\documentclass[11pt]{article}
\usepackage[a4paper,left=2.0cm,right=2.0cm,top=3cm,bottom=3cm]{geometry}
\usepackage{amsmath,amsfonts,amssymb}
\usepackage{mathtools}
\usepackage{siunitx}
\usepackage{fancyhdr}
\usepackage{graphicx}
\usepackage[export]{adjustbox}
\usepackage{float}
\usepackage{subfig}

\usepackage{hyperref}
\hypersetup{hidelinks}

\usepackage{booktabs}
\usepackage{multirow}

\def\@maketitle{\newpage
  \null
  \vskip 1em \begin{center} \let \footnote \thanks
    {\normalsize \bf \@title \par} \vskip 1em {\normalsize \lineskip .5em
	  \begin{tabular}[t]{c} \@author
      \end{tabular}\par}
	\vskip 1em {\normalsize \@date} \end{center}
	\par
  \vskip 1em}

\def\*#1{\mathbf{#1}}
\def\nsd{{n_\mathrm{sd}}}
\def\nx{{\boldsymbol{\nabla}}_{\mathbf{x}}}
\def\nm{{\boldsymbol{\nabla}}_{\#}}

\newcommand{\bd}{\mathbf{d}}

\newcommand{\bn}{\mathbf{n}}

\newcommand{\bq}{\mathbf{q}}

\newcommand{\bu}{\mathbf{u}}

\newcommand{\bw}{\mathbf{w}}
\newcommand{\bx}{\mathbf{x}}

\newcommand{\bA}{\mathbf{A}}

\newcommand{\bH}{\mathbf{H}}
\newcommand{\bI}{\mathbf{I}}

\newcommand{\bK}{\mathbf{K}}

\newcommand{\bU}{\mathbf{U}}

\newcommand{\bW}{\mathbf{W}}

\newcommand{\bY}{\mathbf{Y}}

\newcommand{\zero}{\mathbf{0}}

\newcommand{\strain}{{\boldsymbol{\varepsilon}}}
\newcommand{\stress}{{\boldsymbol{\sigma}}}

\newcommand{\btau}{\boldsymbol{\tau}}
\newcommand{\ddx}[2]{\frac{\partial #1}{\partial #2}}

\def\va3{\left[\begin{array}{c}v_1^a\\v_2^a\\v_3^a\end{array}\right]}
\def\ub3{\left[\begin{array}{c}u_1^b\\u_2^b\\u_3^b\end{array}\right]}
\def\fb3{\left[\begin{array}{c}f_1^b\\f_2^b\\f_3^b\end{array}\right]}

\def\gradv3{\left[\begin{array}{ccc}
	v_{1,x_1}&v_{1,x_2}&v_{1,x_3}\\
	v_{2,x_1}&v_{2,x_2}&v_{2,x_3}\\
	v_{3,x_1}&v_{3,x_2}&v_{3,x_3}\\
\end{array}\right]}
\def\strainu3{\left[\begin{array}{ccc}
	2 u_{1,x_1}&u_{1,x_2}+u_{2,x_1}&u_{1,x_3}+u_{3,x_1}\\
	u_{2,x_1}+u_{1,x_2}&2 u_{2,x_2}&u_{2,x_3}+u_{3,x_2}\\
	u_{3,x_1}+u_{1,x_3}&u_{3,x_2}+u_{2,x_3}&2 u_{3,x_3}
\end{array}\right]}
\def\divstrainv3{\left[\begin{array}{c}
	2v_{1,x_1x_1}+v_{1,x_2x_2}+ v_{1,x_3x_3}+v_{2,x_1x_2}+ v_{3,x_1x_3}\\
	 v_{1,x_2x_1}+v_{2,x_1x_1}+2v_{2,x_2x_2}+v_{2,x_3x_3}+ v_{3,x_2x_3}\\
	 v_{1,x_3x_1}+v_{2,x_3x_2}+ v_{3,x_1x_1}+v_{3,x_2x_2}+2v_{3,x_3x_3}
\end{array}\right]}
\def\divstrainu3{\left[\begin{array}{c}
	2u_{1,x_1x_1}+u_{1,x_2x_2}+ u_{1,x_3x_3}+u_{2,x_1x_2}+ u_{3,x_1x_3}\\
	 u_{1,x_2x_1}+u_{2,x_1x_1}+2u_{2,x_2x_2}+u_{2,x_3x_3}+ u_{3,x_2x_3}\\
	 u_{1,x_3x_1}+u_{2,x_3x_2}+ u_{3,x_1x_1}+u_{3,x_2x_2}+2u_{3,x_3x_3}
\end{array}\right]}
\def\gradq3{\left[\begin{array}{c}q_{x_1}\\q_{x_2}\\q_{x_3}\end{array}\right]}
\def\gradp3{\left[\begin{array}{c}p_{x_1}\\p_{x_2}\\p_{x_3}\end{array}\right]}
\newcommand{\rhou}{\left( \rho \bu \right)}
\newcommand{\rhoe}{\left( \rho e \right)}

\newcommand{\rhoeu}{\left( \rho e \bu \right)}
\newcommand{\drdt}{\frac{\partial \rho}{\partial t}}
\newcommand{\drudt}{\frac{\partial \rhou}{\partial t}}
\newcommand{\dredt}{\frac{\partial \rhoe}{\partial t}}

\let\oldhash\#
\renewcommand{\#}{ {\scalebox{0.5} {\oldhash}} }

\begin{document}

{\large \bf Four-Dimensional Elastically Deformed Simplex Space-Time Meshes for Domains with Time Variant Topology} \par \vspace{6pt}
{\large Max von Danwitz\footnote{Chair for Computational
Analysis of Technical Systems (CATS),
RWTH Aachen University, 52056 Aachen,
Germany. Email: \{danwitz, antony, hosters, behr\}@cats.rwth-aachen.de}, Patrick Antony\footnotemark[1], Fabian Key\footnotemark[1]$^,$\footnote{Institute of Lightweight Design and Structural Biomechanics,
TU Vienna, A-1060 Vienna, Austria. Email: fabian.key@tuwien.ac.at }, Norbert Hosters\footnotemark[1], Marek Behr\footnotemark[1]} \par \vspace{6pt}

\begin{abstract}
\noindent
Thinking of the flow through biological or technical valves, there is a variety of applications in which the topology of a fluid domain changes over time. This topology change is characteristic for the physical behaviour, but poses a particular challenge in computer simulations. A way to overcome this challenge is to consider the space-time extent of the application as a contiguous computational domain. In this work, we obtain a boundary conforming discretization of the space-time domain with four-dimensional simplex elements (pentatopes). To facilitate the construction of pentatope meshes for complex geometries, the widely used elastic mesh update method is extended to four-dimensional meshes. In the resulting workflow, the topology change is elegantly included in the pentatope mesh and does not require any additional treatment during the simulation. The potential of simplex space-time meshes for domains with time variant topology is demonstrated in a valve simulation and a flow simulation inspired by a clamped artery.
\end{abstract}


\section{Introduction and Problem Definition}
\label{sec:intro}
Computer simulations have been widely and successfully used to understand and predict physical behaviour in biological, medical, technical, and many other applications. A central entity in simulation technology is the computational domain on which the solution is sought numerically. Commonly, the computational domain coincides with the spatial extent $\Omega \subset \mathbb{R}^{3}$ of the physical object under consideration. In many applications, this leads to a time-dependent spatial domain $\Omega (t)$ and---thinking for example of valves or bearings---there are also various applications where $\Omega$ changes its topology over time. 

Several approaches have been developed to handle the time variant topology of $\Omega$. Our focus is on mesh-based approaches, which can be broadly organised in three groups, namely, boundary conforming, quasi boundary conforming with residual gap, and non-boundary conforming approaches. Non-boundary conforming (non-interface fitted) approaches find a way to shift the complication of the topology change from the discretisation to the integration. In the field of fluid-structure contact interaction (FSCI) a number of methods has recently been developed. Ager et al.~\cite{ager2019poro} employ a fixed Eulerian mesh as fluid discretization and handle boundary motion and topology changes of the physical fluid domain with an approach based on CutFEM~\cite{burman2015cutfem}. This includes a suitable method for the numerical integration on boundary intersected elements~\cite{sudharkar2014accurate}, a ghost penalty stabilization to overcome issues of very small cuts, and a Nitsche-based weak imposition of fluid boundary (coupling) conditions. Kapada et al. propose a computational framework based on Cartesian hierarchical B-spline grids for the fluid discretization~\cite{kadapa2018stabilised}. Therein, integration over intersected elements is handled via sub-triangulation or uniform subdivision of the rectangular grid; a ghost penalty stabilization and Nitsche's method for fluid boundary conditions is used. Both aforementioned formulations bear similarities with the immersogeometric variational framework proposed previously by Kamensky et al.~\cite{kamensky2015heart}.

The second class of approaches use a boundary conforming discretization, but avoid the actual topology change of the fluid domain with a small residual gap. Ensuring the residual gap with a penalty force or a displacement restriction in the mesh motion, space-time~\cite{sathe2008modeling} or ALE methods~\cite{bogaers2016evaluation} are employed to handle the deforming fluid domain.

Finally, boundary conforming discretizations for spatial computational domains with topology changes can be obtained by reformulating the problem as space-time problem with contiguous computational domain. For two-dimensional spatial computational domains, this approach has been implemented with unstructured finite volume~\cite{rendall2012conservative} and finite element meshes~\cite{danwitz2019simplex}. 
 
\begin{figure}
\centering
\includegraphics{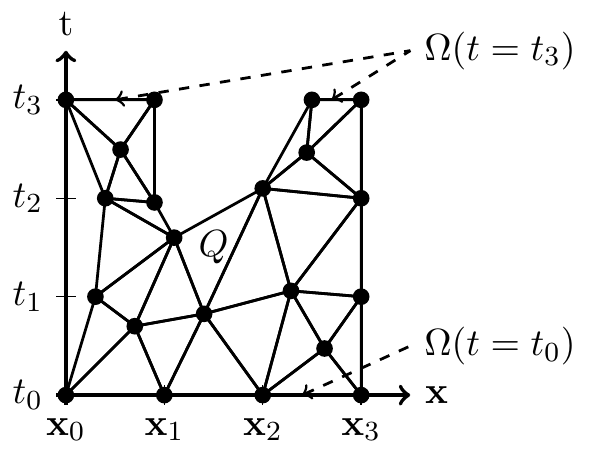}
\caption{Simplex space-time (SST) mesh of spatial computational domain with time variant topology.}
\label{fig:ST}
\end{figure}

In this work, we will follow the boundary conforming approach and choose the space-time domain $Q \subset \mathbb{R}^4$ as computational domain. The concept is visualised in Figure~\ref{fig:ST} with the spatial domain collapsed onto the horizontal axis. Between $t_1$ and $t_2$, the spatial domain $\Omega$ splits into two unconnected parts, however, the computational space-time domain $Q$ is contiguous. A boundary conforming discretization with pentatope elements leads to a simplex space-time (SST) mesh, which is fully unstructured in space and time.

To generate such four-dimensional simplex meshes, several similar approaches have been developed recently. The robust meshing strategy proposed by Behr~\cite{behr2008simplex} extrudes a tetrahedral mesh resulting in a four-dimensional mesh of tensor product elements. These hyperprisms are split into pentatopes with an element-wise Delauney triangulation of the perturbed nodal coordinates. In a further development, Karabelas and Neum\"{u}ller replace the element-wise Delauney triangulation with a predefined decomposition of each hyperprism requiring a consistently numbered tetrahedral mesh~\cite{karabelas2019generating}.
A third approach of Wang similarly employs a global node indexing scheme and extends it with a node insertion procedure to support local mesh density operations~\cite{wang2015discontinuous}. The above mentioned strategies have in common that the four-dimensional mesh is based on an extruded tetrahedral mesh. An alternative approach to generate high-quality pentatope meshes could be based on Coxeter triangulations, with the currently open issue of generating boundary conforming meshes~\cite{choudhary2020coxeter}.

In most cases, SST meshes are employed to facilitate adaptive mesh refinement in space and time. Suitable pentatope mesh refinement procedures have been explored by Neum\"{u}ller and Steinbach~\cite{neumuller2011refinement} and Grande~\cite{grande2019red}. Anisotropic four-dimensional mesh adaptation is pioneered by Caplan et al~\cite{caplan2017anisotropic,caplan2020four} and successfully employed in the solution of the advection diffusion equation~\cite{caplan2019four}. Further, recent application examples of four-dimensional SST meshes from the field of mathematics deal with parabolic evolution problems~\cite{langer2019coeff,steinbach2019space} or a broader class of transient PDEs recast as constrained first-order system~\cite{voronin2018space}. In the field of computational engineering, adaptive temporal refinement of pentatope meshes is used for two-phase flow simulations~\cite{karyofylli2018simplex}---also combined with complex material laws such as the Carreau-Yasuda-WLF model~\cite{karyofylli2019simplex} or the $\mu(I)$-rheology~\cite{gesenhues2020advanced}---as well as gas flow simulations in the piston ring pack of internal combustion engines~\cite{danwitz2019simplex}.

In this work, we generate pentatope finite element meshes for spatial domains with time variant topology. Therefore, we combine the extrusion based approach by Behr~\cite{behr2008simplex} with a four-dimensional extension of the elastic mesh update method (EMUM). Originally proposed as elastic grid approach by Lynch~\cite{lynch1982unified}, EMUM was refined and employed as automatic mesh moving scheme for the deforming spatial domain / stabilized space-time finite element formulation~\cite{johnson1994mesh}. Since then, it has been widely used and become a standard technique to handle moving domains in fluid-structure interaction (FSI) simulations~\cite{hubner2004monolithic,bazilevs2008isogeometric,bazilevs2013computational}. Also, very recent FSI simulations rely on it~\cite{spenke2020multi,la2020role,liu2020fluid}. Furthermore, EMUM has been used in the context of free-surface flows (for example in~\cite{zwicke2017boundary}) and to update finite element meshes according to prescribed boundary displacements (for example in~\cite{wendling2019cfd}). The four-dimensional extension of EMUM (4DEMUM) presented in Section~\ref{sec:method} allows us to obtain boundary conforming pentatope meshes of complex geometries. The geometry may have holes, and does not have to have a tensor product shape in any dimension. This means that the pentatope mesh can account for a time variant topology of the spatial domain. However, the mesh generation method is still limited to cases where the four-dimensional geometry can be obtained by extrusion of a (complex) three-dimensional geometry with subsequent elastic deformation (in the sense of 4DEMUM).

Now, we proceed as follows. In Section~\ref{sec:method}, we describe the strong and weak form for a finite element implementation of 4DEMUM. Next, Section~\ref{sec:workflow} presents the embedding of 4DEMUM into the workflow of a four-dimensional space-time finite element simulation, followed by Section~\ref{sec:examples} which aims at experimentally validating the method, as well as demonstrating the particular potential of simulations on SST meshes for domains with time variant topology. Finally, Section~\ref{sec:conclusion} contains concluding remarks and a brief outlook.

\section{Four-Dimensional Elastic Mesh Update Method}
\label{sec:method}
\begin{figure}[h]
     \centering
     \subfloat[Mesh on extruded domain\label{fig:QMesh}]{%
      \hspace{1.5cm}
	\includegraphics{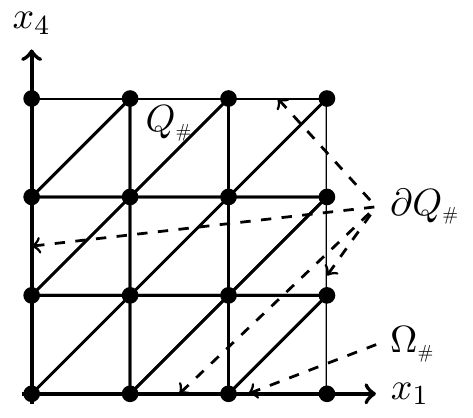}
  \hspace{1.5cm}} \quad 
 \subfloat[Updated mesh on physical space-time domain\label{fig:QPhys}]{%
  \hspace{1.5cm}
	\includegraphics{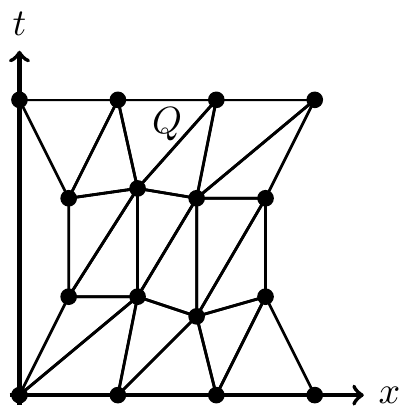}
         \hspace{1.5cm} }
     \caption{Four-dimensional elastic mesh update method.}
     \label{fig:E4DMUM}
 \end{figure}
 
The basic concept of the elastic mesh update method is to update the node positions according to the deformations of an elastic solid. As outlined above, the method has been widely used to adjust two- or three-dimensional meshes to deforming boundaries. In the following, we extend the method to four-dimensional meshes. Therefore, we consider a virtual four-dimensional linear elastic solid occupying the region $Q_\# \subset \mathbb{R}^4$. 

Figure~\ref{fig:QMesh} displays two dimensions of an example of $Q_\#$. Let's assume the domain boundary $\partial Q_\#$ consists of a Neumann part $\partial Q_\#^N$ and a Dirichlet part $\partial Q_\#^D$. For each degree of freedom the boundary is split such that  $\partial Q_\# = \partial Q_\#^D \cup \partial Q_\#^N$  and $ \partial Q_\#^D \cap \partial Q_\#^N = \emptyset $. We can now formulate the four-dimensional elastostatic problem
\begin{equation}
	\nm \cdot \stress( \bd ) = \zero, \quad \mbox{on $Q_\#$}, \quad   \bd = \bd_D \quad \mbox{on $\partial Q_\#^D$},
    \label{eq:e-linear}
\end{equation}
with Dirichlet boundary conditions $\bd_D$ prescribed on $\partial Q_\#^D$. Where no Dirichlet values are prescribed, a homogeneous Neumann boundary condition is assumed.

The assumption of a homogeneous, isotropic, linear elastic material behaviour leads to a constitutive equation 
\begin{equation}
\stress = \lambda \left[ \hbox{tr}\left(\strain\right) \right] \bI + 2 \mu \strain,
\end{equation}  
where $\stress$ is an equivalent of the Cauchy stress tensor and $\lambda$ and $\mu$ are the Lam\'{e} parameters. In this work, we choose $\lambda = \mu = 1$. Further, $\bI \in \mathbb{R}^{4 \times 4}$ denotes the identity matrix and $\hbox{tr} (\strain)$ is the sum of the main diagonal entries of the strain tensor.

The solid's deformation is measured with a linear strain model
\begin{equation}
\strain = \frac{1}{2} \left[ \nm \bd + \left( \nm \bd \right)^\top \right],
\end{equation}
where $\bd (\bx_{\#}) \in \mathbb{R}^4$ are displacements in the four spatial directions. The displacements are observed at the nodes of the extruded mesh and the nodal positions are described by their coordinates $\bx_{\#} = (x_1, x_2, x_3, x_4)^T \in Q_\#$. Accordingly, the nabla operator $\nm = \left[\ddx{(.)}{x_1},\ddx{(.)}{x_2},\ddx{(.)}{x_3},\ddx{(.)}{x_4}\right]^T$ collects the partial derivatives. 

For an admissible test function space $V_{h}$ with functions that vanish on $\partial Q_\#^D$, and a suitable trial function space $S_h$ that respects the Dirichlet boundary conditions, the discretised weak form of Equation~\eqref{eq:e-linear} can be stated as: 
Find $\bd^h \in S_{h}$, such that for all $\bw^h \in V_{h} $
\begin{align}
\label{eq:weakMesh}
0 = \int_{Q} \nm \bw^h \,  : \, \stress(\bd^h) \,dQ  = \int_{Q} \nm \bw^h \,  : \left\{ \lambda \hbox{tr} \left( \nm \bd^h \right) \bI + \mu \left[ \nm \bd^h + \left( \nm \bd^h \right)^\top \right] \right\} dQ.
\end{align}
For details on pentatope Lagrange finite elements the reader is referred to the previous publications~\cite{danwitz2019simplex,behr2008simplex}. In short, the weak form above results in a linear equation system that is solved for the displacements $\bd^h$. Finally, the node positions with coordinates $\bx = (x, y, z, t)^T \in Q$ are obtained by adding the displacements $\bd^h$ to the extruded mesh coordinates:
\begin{equation}
\bx = \bx_{\#} + \bd^h.
\end{equation}
In this addition, we identify $(x,y,z,t)$ with $(x_1, x_2, x_3, x_4)$. Two dimensions of a resulting mesh on the space-time domain are shown in Figure~\ref{fig:QPhys}. The displacements in all four dimensions are coupled, so interior mesh nodes can move in any direction, even if non-zero boundary displacements are applied only to specific ones.

\section{Simulation Workflow}
\label{sec:workflow}
\begin{figure}
\centering
\includegraphics[width=\textwidth]{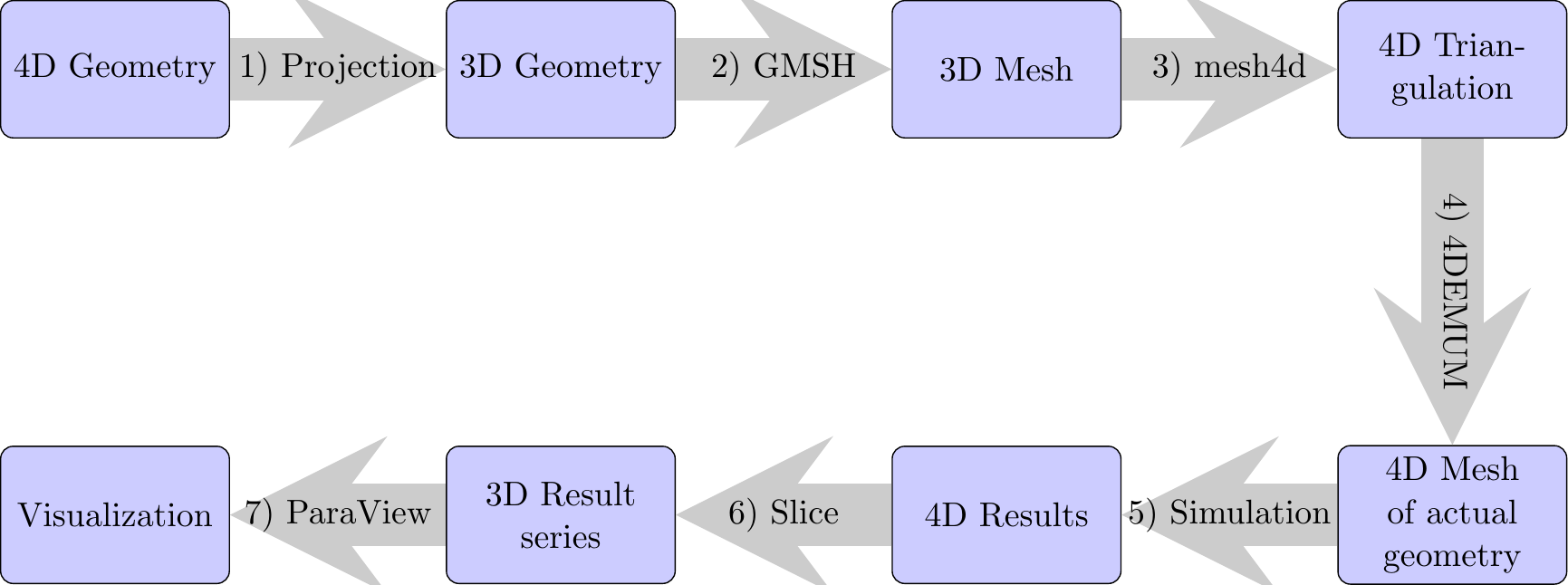}
\caption{Workflow for four-dimensional SST finite element simulation.}
\label{fig:workflow}
\end{figure}

An overview of the workflow of a four-dimensional SST finite element simulation is given in Figure~\ref{fig:workflow}. The arrows correspond to seven steps that are taken to complete a specific subtask in the workflow. When suitable, the arrows are labeled with the software employed to complete the task. The first four steps are considered pre-processing and produce the SST meshes covering the considered physical space-time domain $Q$. In detail, they are: First, identify a projection $\Omega_\# = \mathbb{P} (Q) \subset \mathbb{R}^3$, such that complex features of $Q$ are included in $\Omega_\#$. In our work, this step is performed manually. For an example of such a projection result see Figure~\ref{fig:arteryMesh}.  As second step, generate a tetrahedral mesh covering $\Omega_\#$, for example with GMSH~\cite{geuzaine2009gmsh}. Then, according to Behr~\cite{behr2008simplex}, extrude the mesh covering $\Omega_\#$ over an interval $I \subset \mathbb{R}$ resulting in a mesh of tensor-product elements and apply an element-wise Delauney triangulation to generate a pentatope mesh covering $Q_\# = \Omega_\# \times I $ (third step). The fourth steps is to use 4DEMUM as described in Section~\ref{sec:method} to deform the pentatope mesh to cover $Q$. Note that the time dimension can be one of the initial three dimensions, and extrusion can be used to generate one of the spatial dimensions of the mesh. 

With a boundary conforming discretization at hand, the fifth step is to perform the space-time simulation. Details of the simulation depend on the application example and are given in Section~\ref{sec:examples}. In any case, the solution of the finite problem is obtained on the unstructured space-time mesh. To visualise the results over time, a series of tetrahedral meshes covering the spatial domain at given time instances is generated. The node positions of these meshes are passed as query points to an efficient interpolation tool~\cite{fernandez2019mesh}. The tool identifies the pentatope of the space-time mesh which contains the query point  and performs a linear barycentric interpolation on the element (step 6). This interpolation is in accordance with the finite element approximation of the simulation. For domains with convex boundaries some of the query points can be located slightly outside of the space-time mesh. In such cases, the solution from the element with the closest center is extrapolated. As final step 7, the data on the tetrahedral meshes can be easily visualised with available tools such as ParaView~\cite{ayachit2015paraview}. An alternative visualisation approach described by Karabelas and Neum\"{u}ller~\cite{karabelas2019generating} is based on the element-wise intersection of the pentatope mesh with a hyperplane.

\section{Application Examples}
\label{sec:examples}
For the physical simulations in the following application examples, we distinguish between space and time coordinates. The spatial nabla operator $\nx  =  \ddx{(\cdot)}{\bx}$ collects the derivatives with respect to the spatial coordinates. In the following two examples, we consider the compressible Navier-Stokes equations
\begin{eqnarray}
	\label{eqn-cont}
	\drdt + \nx \cdot \rhou = 0
	\quad &\mbox{on $Q,$} \\
	\label{eqn-mom}
	\drudt +\nx \cdot \left[\rhou \otimes \bu\right] + \nx \, p - \nx \cdot \btau = \zero
	\quad &\mbox{on $Q,$} \\
	\label{eqn-energy}
	\dredt +\nx \cdot \rhoeu + \nx \cdot (p \bu) 
	- \nx \cdot (\btau \bu) + \nx \cdot  \bq = 0
	\quad &\mbox{on $Q.$} 
\end{eqnarray}
Here $\rho(\bx,t)$, $\bu(\bx,t)$, $p(\bx,t)$, $\btau(\bx,t)$, $e(\bx,t)$,
and $\bq(\bx,t)$ are density, 
velocity vector, pressure, viscous stress tensor, total energy per unit mass, 
and heat flux vector, respectively. For convenience, the conservation variables are collected in the vector $\mathbf{U} = \left[  \rho, \rho \bu, \rho e \right]^\top$. We consider an ideal, calorically perfect gas, with the specific gas constant $R=\SI{287}{\joule \per \kilo \gram \per \kelvin}$, the ratio of specific heats $\gamma = 1.4$,  and a Prandtl number of $\operatorname{Pr} = \frac{\nu}{\kappa}\frac{\gamma R}{\gamma-1} = 0.71 $. With these assumptions, we can perform a change of variables from the conservation variables $\bU = \bU(\bY)$ to the pressure-primitive variables $\mathbf{Y} = \left[  p,  \bu, T \right]^\top$. The latter are used as primary unknowns and the governing equations above are formulated as generalized advective-diffusive system
\begin{equation}
\textbf{Res}(\bY) \coloneqq \bA_0 \bY_{,t} + \left( \bA^{\text{adv\textbackslash p}}_i + \bA^{\text{p}}_i + \bA^{\text{sp}}_i \right) \bY_{,i} - (\bK_{ij} \bY_{,j})_{,i} = \zero, \quad i,j = 1, \dots , \nsd.
\label{eq:genad}
\end{equation}
Therein, partial time derivatives are denoted with $(\cdot)_{,t}$, partial derivatives in each of the $\nsd = 3$ spatial directions are denoted with $(\cdot)_{,i}$, and Einstein summation convention applies to repeated indices. Details of the generalized advection matrices $\bA$ and diffusion matrices $\bK$ can be found in~\cite{danwitz2019simplex, xu2017compressible}.

\begin{figure}
\centering
\includegraphics{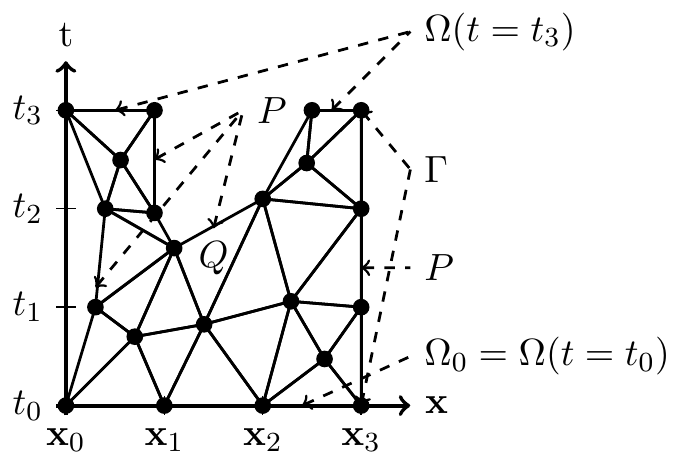}
\caption{Physical space-time domain and labeled boundary parts.}
\label{fig:PhysicalST}
\end{figure}

To find a numerical solution to the compressible Navier-Stokes equations with space-time finite elements, a discretisation of the physical space-time domain $Q \subset \mathbb{R}^{\nsd+1}$ is required. An example of such a discretisation can be seen in Figure~\ref{fig:PhysicalST}. In practise, a procedure is required to numerically integrate over $Q$, over the spatial computational domain $\Omega_0 \subset \mathbb{R}^{\nsd}$ at the initial time $t_0$, and over $P$ which is the temporal evolution of the spatial domain boundary $\Gamma \subset \mathbb{R}^{\nsd-1}$. For each degree of freedom, the space-time boundary $P$ consist of non-overlapping Dirichlet $P_D$ and Neumann parts $P_N$.

Based on the domains introduced above and for admissible test and trial function spaces $V_{h}$ and $S_h$, the weak form of Equation~\eqref{eq:genad} can be stated as: 
For given initial conditions $\bY^h_0$, find $\bY^h \in S_{h}$, such that for all $\bW^h \in V_{h} $
\begin{align}
\label{eq:weak}
0 & = \int_{Q} \bW^h \cdot \left[  \bA_0 \bY^h_{,t}  +  \left( \bA^{\text{adv\textbackslash p}}_i + \bA^{\text{sp}}_i \right) \bY^h_{,i} \right]  dQ \\ \nonumber
   & + \int_{Q} \bW^h_{,i} \cdot \left[ \bK_{ij} \bY^h_{,j} - \bA_i^{\text{p}} \, \bY^h \right] dQ \\ \nonumber
   & -  \int_{P_N} \bW^h \cdot \bH(\bY^h) \, dP\\ \nonumber
   & + \int_{\Omega_{0}} \bW^h \cdot \left[ \bU(\bY^h) - \bU(\bY^h_0) \right] d\Omega \\ \nonumber
   & + \int_{Q} \!	\left[ (\hat{\bA}_{m})^{T} \bW^h_{,m}	\right] \cdot \btau_{\mathrm{SUPG}} \, \mathbf{Res(Y}^h)\, dQ, \\ \nonumber
   & i,j = 1, \dots , \nsd, \quad m =1, \dots, \nsd+1.
\end{align}
The first integral collects terms of the residual that are multiplied with the test function vector and integrated over $Q$. In the second integral, the derivatives with respect to the spatial coordinates are shifted to the test function vector, leading to the third integral over the Neumann part of the space-time boundary $P_N$. Therein, the boundary normal fluxes $\bH$ are evaluated with the outwards pointing surface normal $\bn = (n_1, n_2, n_3)^T $ as
\begin{equation}
	     \bH = \left[  \begin{array}{c}
	     0 \\
	     -p n_1 +\tau_{1i} n_i\\
	     -p n_2 +\tau_{2i} n_i \\
	     -p n_3 + \tau_{3i} n_i \\
	    -q_i n_i
	     \end{array} \right].
\end{equation}
The initial condition $\bY^h_0$ is weakly enforced with the integral over $\Omega_{0}$. The weak form is completed with a SUPG operator to overcome instabilities of the pure Galerkin formulation which occur in convection-dominated flow simulations. For more details on the space-time finite element scheme for the compressible Navier-Stokes equations, the reader is referred to~\cite{danwitz2019simplex}. 

\subsection{Transient Gas Flow through Valve}
\label{ssec:valve}
\begin{figure}
\subfloat[Side view on fluid domain in open configuration and gap height during valve cycle\label{fig:valveSide}]{
\includegraphics{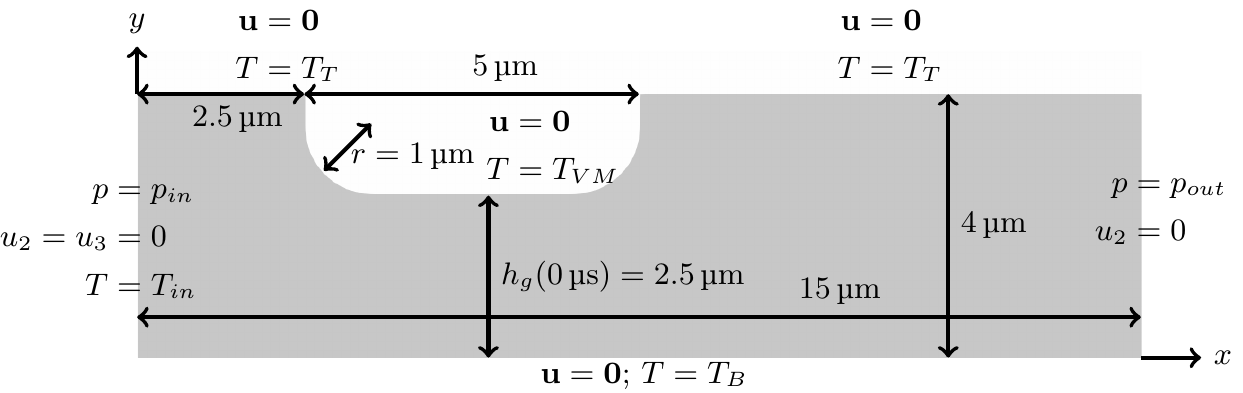} \includegraphics{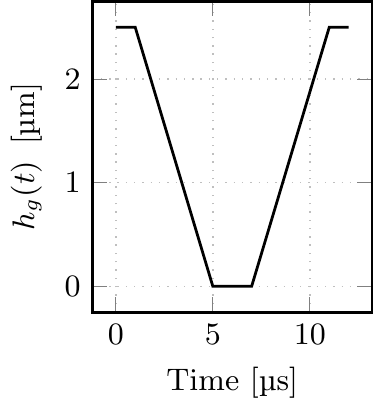}}
\\
\subfloat[Top view on fluid domain with variable exit width $h_e(t)$\label{fig:valveTop}]{
\includegraphics{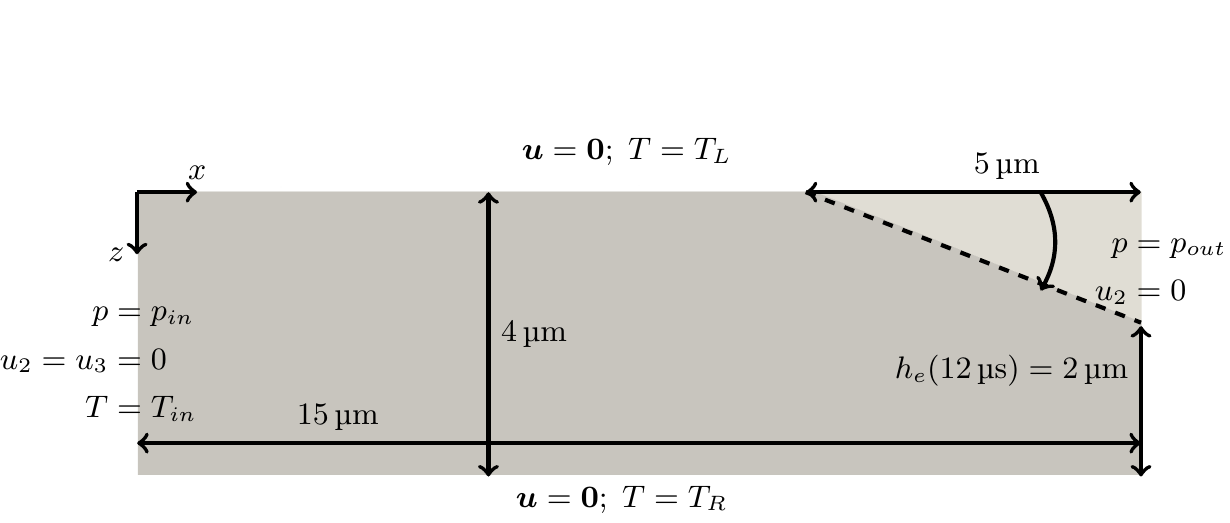} \includegraphics{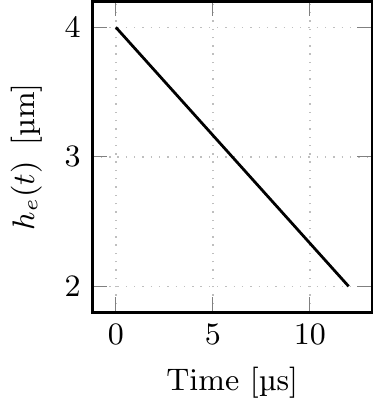}}
\caption{Valve test case. Fluid domain geometry specifications and boundary conditions.}
\label{fig:valveSpec}
\end{figure}

In the following, we present a transient three-dimensional gas flow simulation through a valve. The fluid domain consist of a square channel with \SI{4}{\micro\meter} sides in the $y$-$z$-plane and a length of \SI{15}{\micro\meter} in $x$-direction. In the course of the simulation, the position of the rounded valve member determines the topology of the spatial computational domain (see Figure~\ref{fig:valve4D}). During the valve cycle, the valve member is first lowered and later lifted again. This leads to a splitting of the fluid domain at  \SI{5}{\micro\second} and a reconnection at \SI{7}{\micro\second} (see Figure~\ref{fig:valveSide}). Additionally, the exit cross-section of the channel deforms from the initial square geometry into a rectangular cross-section with half the size. Detailed geometry specifications and boundary conditions are collected in Figure~\ref{fig:valveSpec}. Along the solid walls, no-slip boundary conditions are enforced, i.e., the velocity is set to zero and a wall temperature is prescribed. On the left open boundary, the pressure and temperature are given; on the right open boundary a pressure value is set. Prescribed temperature and pressure values are summarised in Table~\ref{tab:valveBc}. The gas viscosity is modelled using Sutherland's relation
\begin{equation}
\nu = \nu_{ref} \frac{T_{ref}+C}{T+C}\left(\frac{T}{T_{ref}}\right)^{\frac{3}{2}},
\label{eq:Sutherland}
\end{equation}
with $\nu_{ref} =\SI{21.7e-6}{\pascal\second}$, $T_{ref}=  373.15$ K, $C=120$ K.

\newcommand{\myL}{0.6\textwidth}

\begin{table}
\centering
\caption{Valve test case. Boundary conditions.}
\begin{tabular}{c c c c c c c}
\toprule
 $T_{in} \, [\SI{}{\kelvin}]$ & $T_{T} \, [\SI{}{\kelvin}]$ & $T_{VM} \, [\SI{}{\kelvin}]$ & $T_{B} \, [\SI{}{\kelvin}]$ & $T_{L}=T_{R} \, [\SI{}{\kelvin}]$ & $p_{in} \, [\SI{}{\pascal}]$  & $p_{out} \, [\SI{}{\pascal}]$ \\
  \midrule
373.15 & 383.15 & 403.15 & 393.15 & $T_B - \frac{10.0 y}{\SI{4}{\micro\meter}}$ & \SI{1.3e5}{}  & \SI{1.0e5}{} \\
 \bottomrule
\label{tab:valveBc}
\end{tabular}
\end{table}

Following the workflow outlined in Section~\ref{sec:workflow}, we start with a projection of four-dimensional space-time geometry of the valve into the $x$-$y$-$t$-hyperplane. The resulting three-dimensional $x$-$y$-$t$-domain corresponds to the temporal evolution of the side view shown in Figure~\ref{fig:valveSide} and includes the topology change of the flow domain. We discretise the $x$-$y$-$t$-domain with \SI{177 708}{} tetrahedral elements. This is sufficient to resolve a transient two-dimensional flow field as documented in~\cite{danwitz2019simplex}. In the next step, the mesh is extruded in z-direction and triangulated with \SI{14 216 640}{} pentatopes, such that both sides of the square channel cross-section are split into 20 line elements. To account for the closing exit, we identify $(x,y,z,t)$ with $(x_1, x_2, x_3, x_4)$ and apply 4DEMUM with 
\begin{equation}
    \bd_D =  \left[  \begin{array}{c}
	     0 \\
	     0 \\
	      \SI{2}{\micro \meter} \\
	     0  \end{array} \right] \cdot%
	     \mathcal{H} \left( x_1-\SI{10}{\micro \meter}\right) \left( \frac{x_1}{\SI{5}{\micro \meter}} - 2 \right) \cdot%
	    \left (1 - \frac{x_3}{\SI{4}{\micro \meter}} \right)\cdot%
	    \frac{x_4}{\SI{12}{\micro \meter}}
    \label{eq:valveMeshBc}
\end{equation}
prescribed as Dirichlet boundary condition on the domain boundary. Therein, $\mathcal{H}(x)$ denotes the Heaviside function. After the mesh update, the fourth dimension of the resulting mesh is interpreted as time and the SST flow simulation is performed. The simulation took 30 minutes of wall clock time on 240 cores using a distributed memory parallelisation based on MPI.

The results displayed in Figure~\ref{fig:valve4D} show the pressure, velocity, and temperature distribution in the closed, half open, and fully opened valve. In Figure~\ref{fig:valve3V12}, it can be clearly seen that the reduced cross-section of the closing exit accelerates the flow to the maximum velocity of \SI[per-mode=symbol]{94}{\meter\per\second}. The influence of the flow on the temperature distribution is shown in Figure~\ref{fig:valve3T12}. In summary, we observe a laminar flow with a maximum Reynolds number $\operatorname{Re} = \frac{\rho\, u_{max}\, h_e(\SI{12}{\micro\second})}{\mu_{ref}} \approx 10$ and Mach number $\operatorname{Ma} = \frac{u_{max}}{\sqrt{\gamma R T}} \approx 0.24$ at $t=\SI{12}{\micro \second}$. The density variations between \SI[per-mode=symbol]{0.89}{\kilo\gram\per\meter\tothe{3}} and \SI[per-mode=symbol]{1.2}{\kilo\gram\per\meter\tothe{3}} underline the necessity to consider compressibility effects. 

\newcommand{\myW}{0.325\textwidth} 
\begin{figure}
\centering
\includegraphics[width=\myW,trim={17cm 9cm 5cm 40cm},clip]{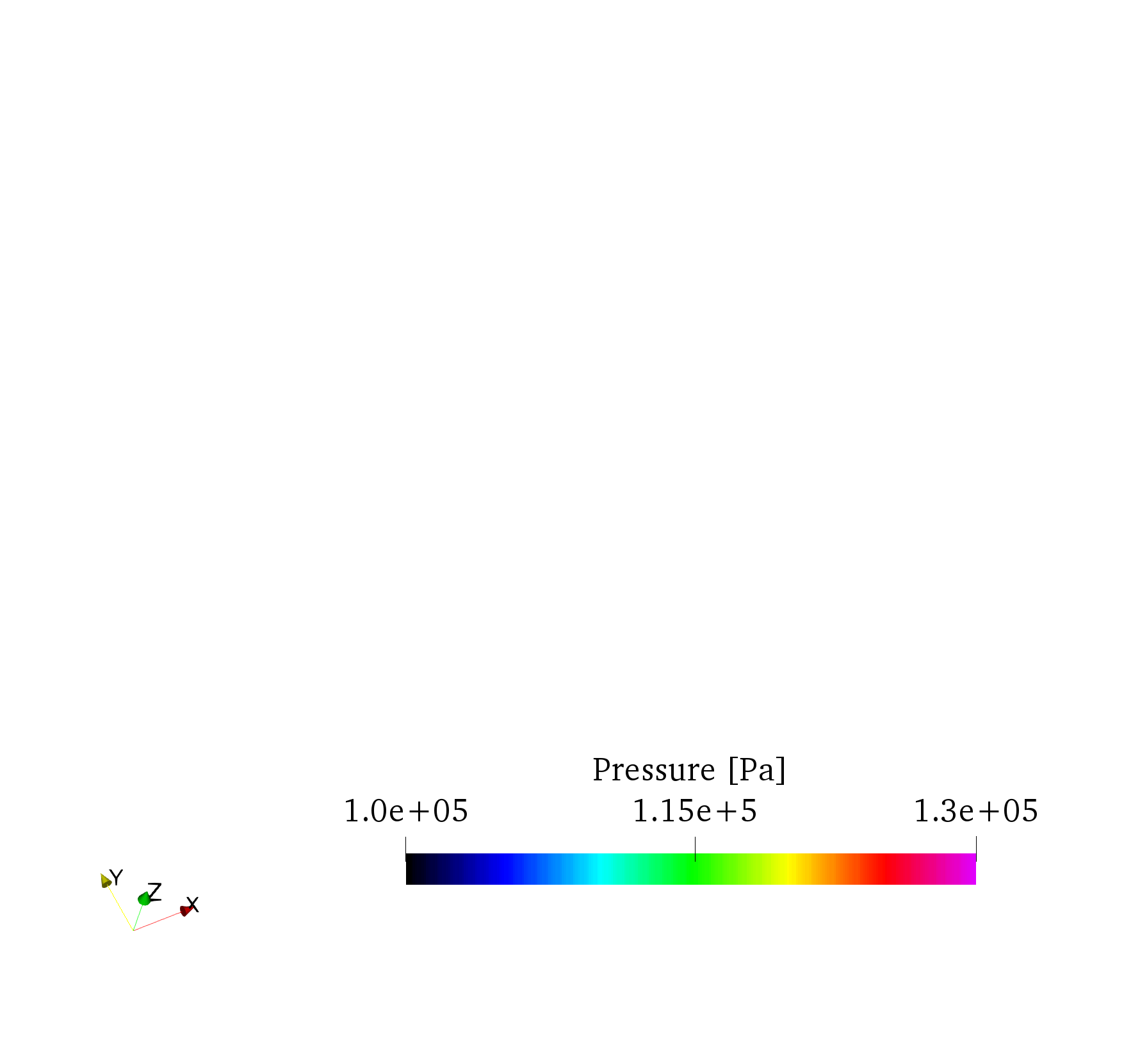}
\includegraphics[width=\myW,trim={17cm 9cm 5cm 40cm},clip]{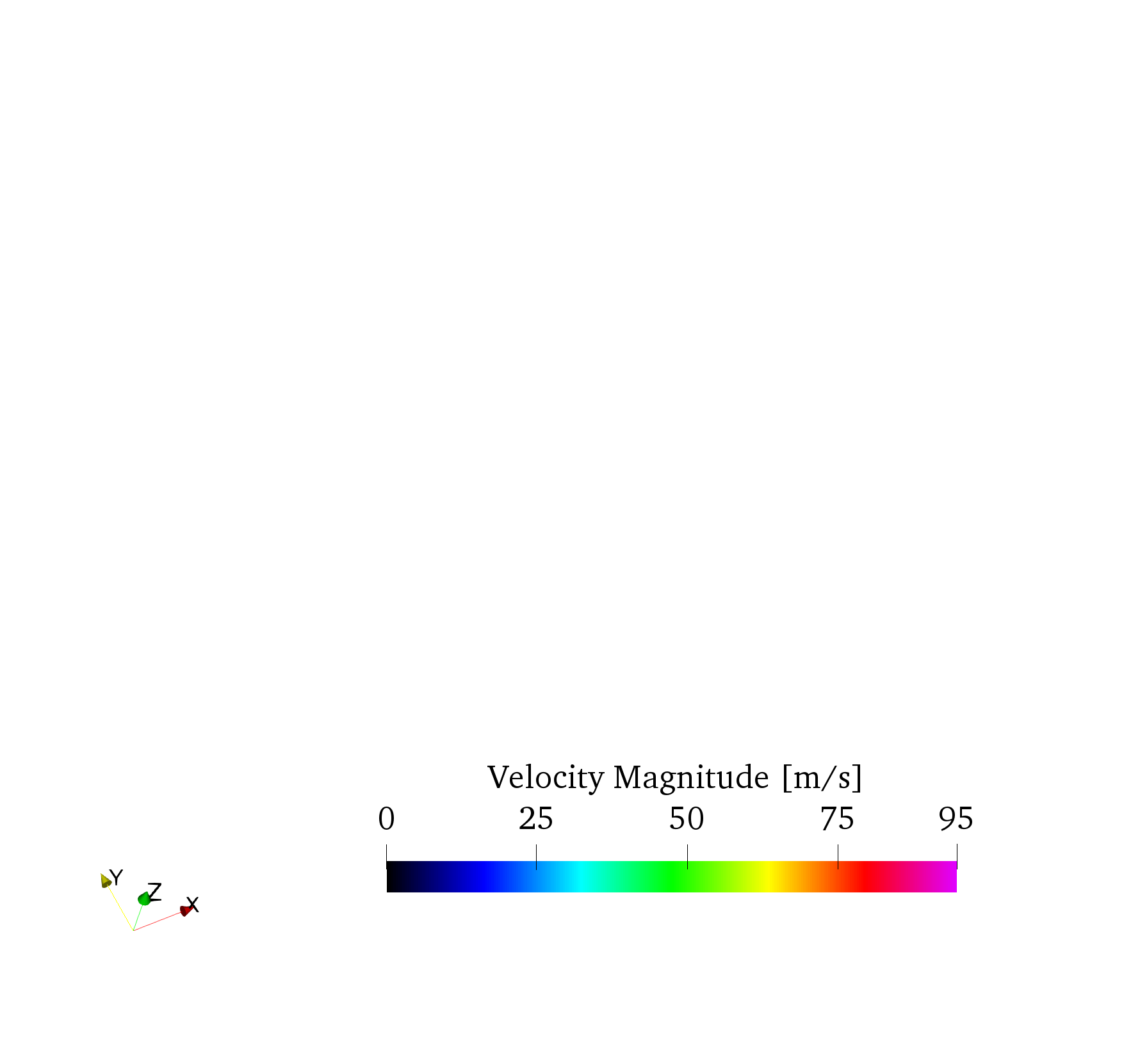}
\includegraphics[width=\myW,trim={17cm 9cm 5cm 40cm},clip]{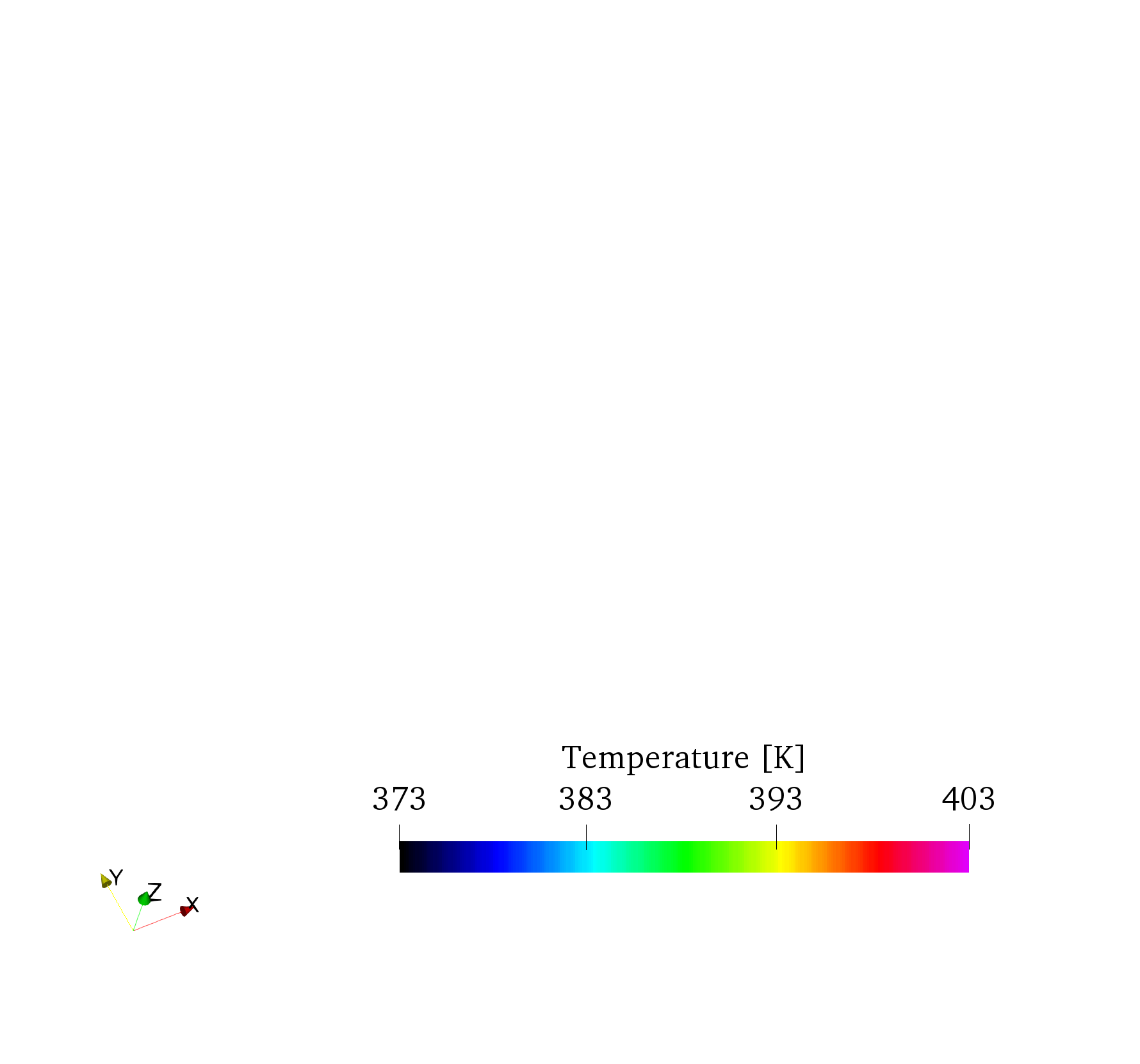}

\subfloat[Pressure distribution at $t=6\, \mu$s \label{fig:valve3P6}]{\includegraphics{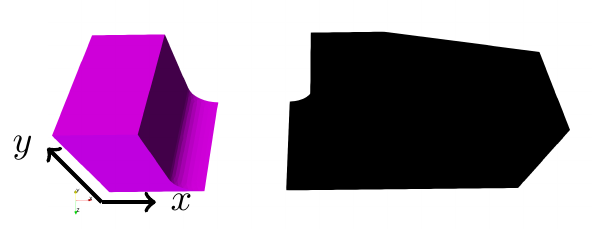}}
\subfloat[Velocity field at $t=6\, \mu$s \label{fig:valve3V6}]{\includegraphics[width=\myW,trim={0cm 0cm 0cm 0cm},clip]{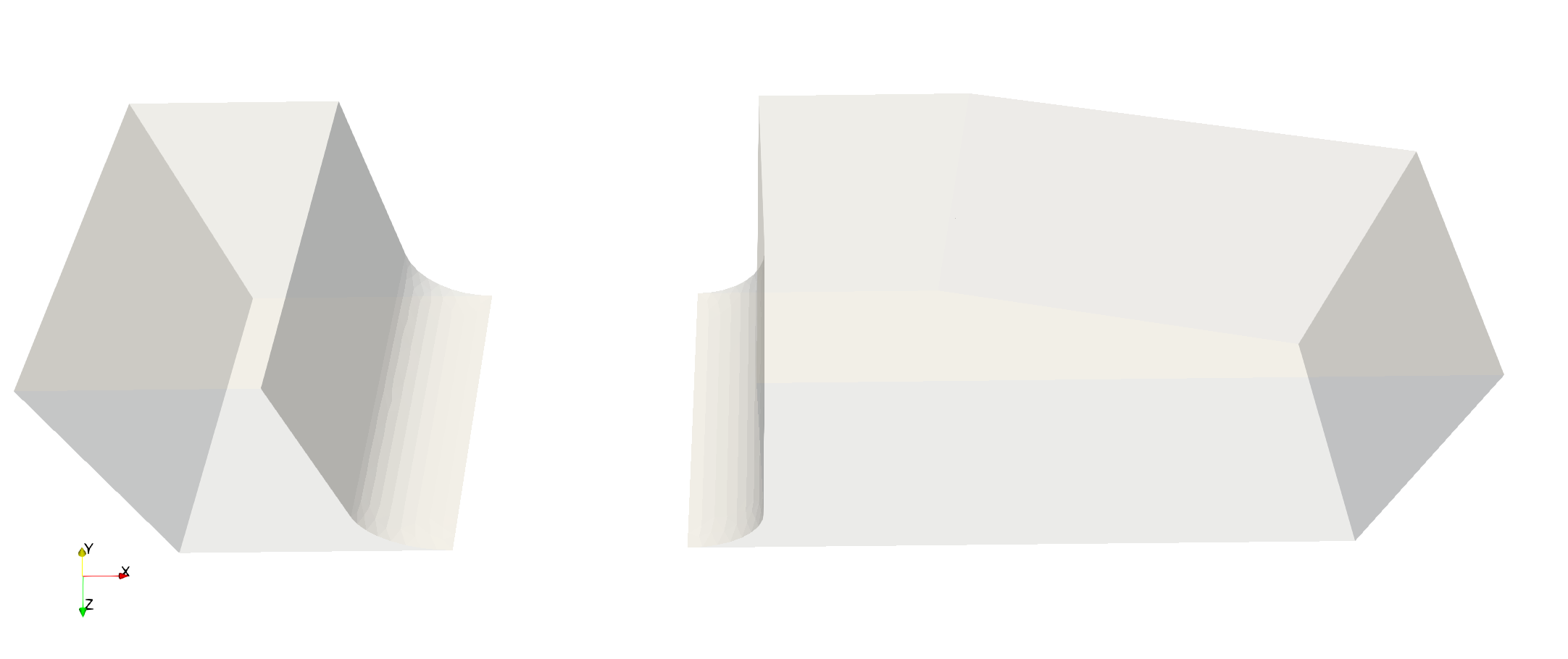}}
\subfloat[Temperature at $z = 3 \, \mu$m, $t=6\, \mu$s \label{fig:valve3T6}]{\includegraphics[width=\myW,trim={0cm 0cm 0cm 0cm},clip]{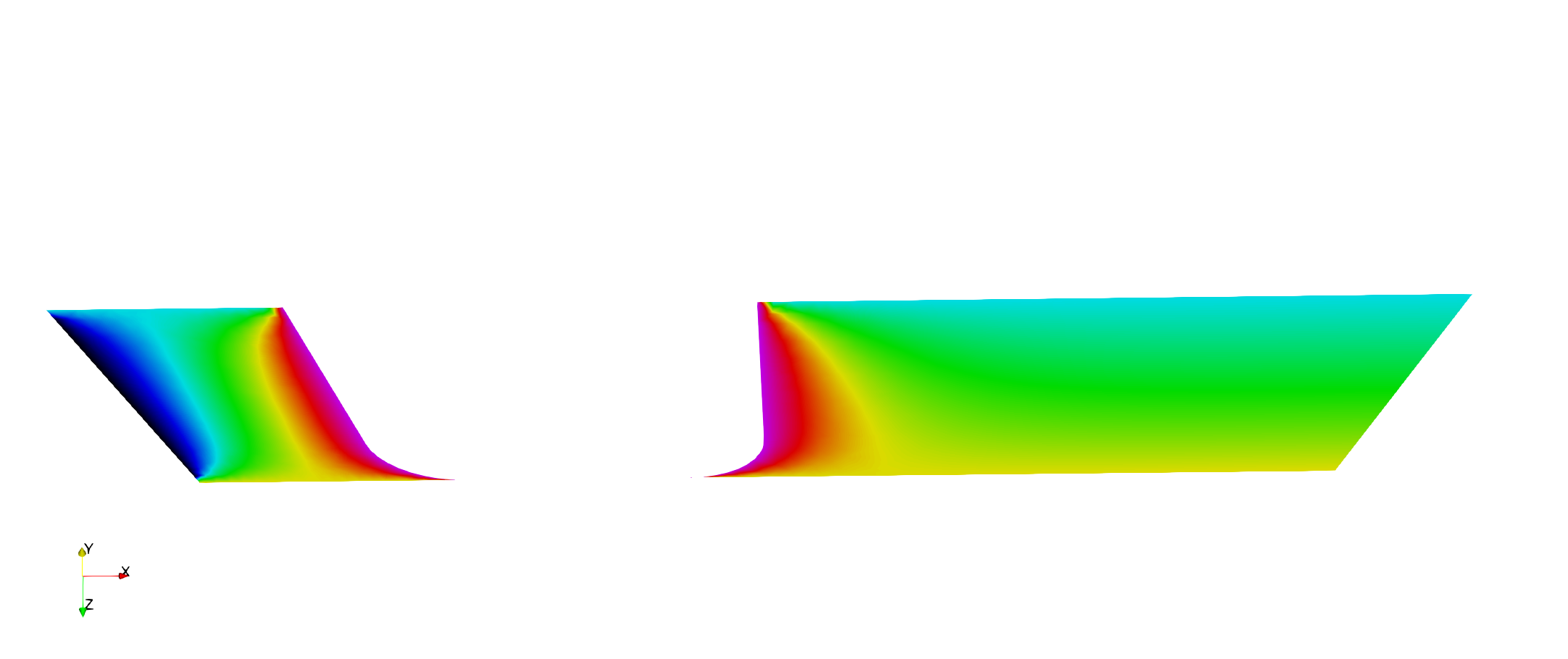}}

\subfloat[Pressure distribution at $t=9\, \mu$s \label{fig:valve3P9}]{\includegraphics{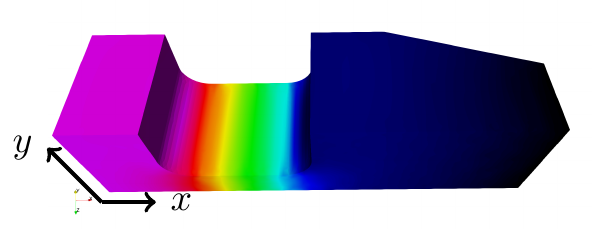}}
\subfloat[Velocity field at $t=9\, \mu$s \label{fig:valve3V9}]{\includegraphics[width=\myW,trim={0cm 0cm 0cm 0cm},clip]{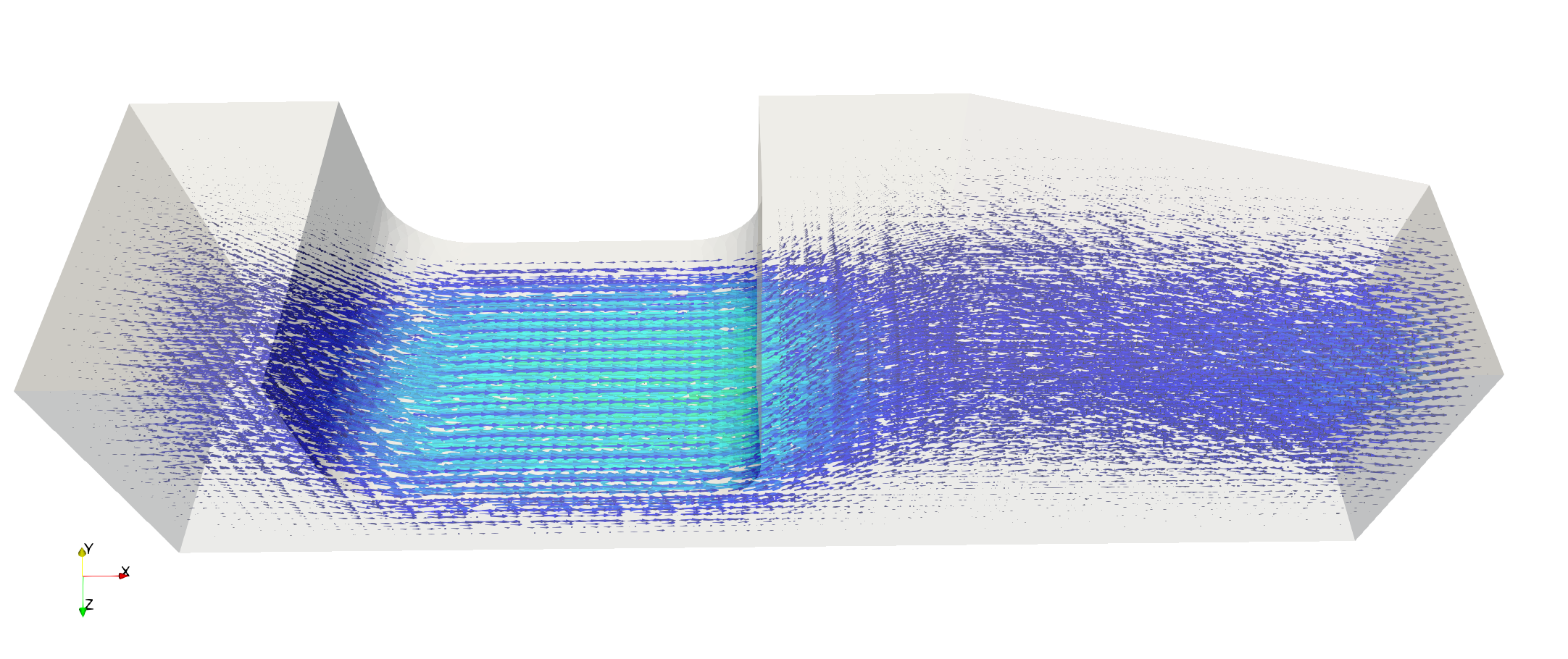}}
\subfloat[Temperature at $z = 3 \, \mu$m,  $t=9\, \mu$s \label{fig:valve3T9}]{\includegraphics[width=\myW,trim={0cm 0cm 0cm 0cm},clip]{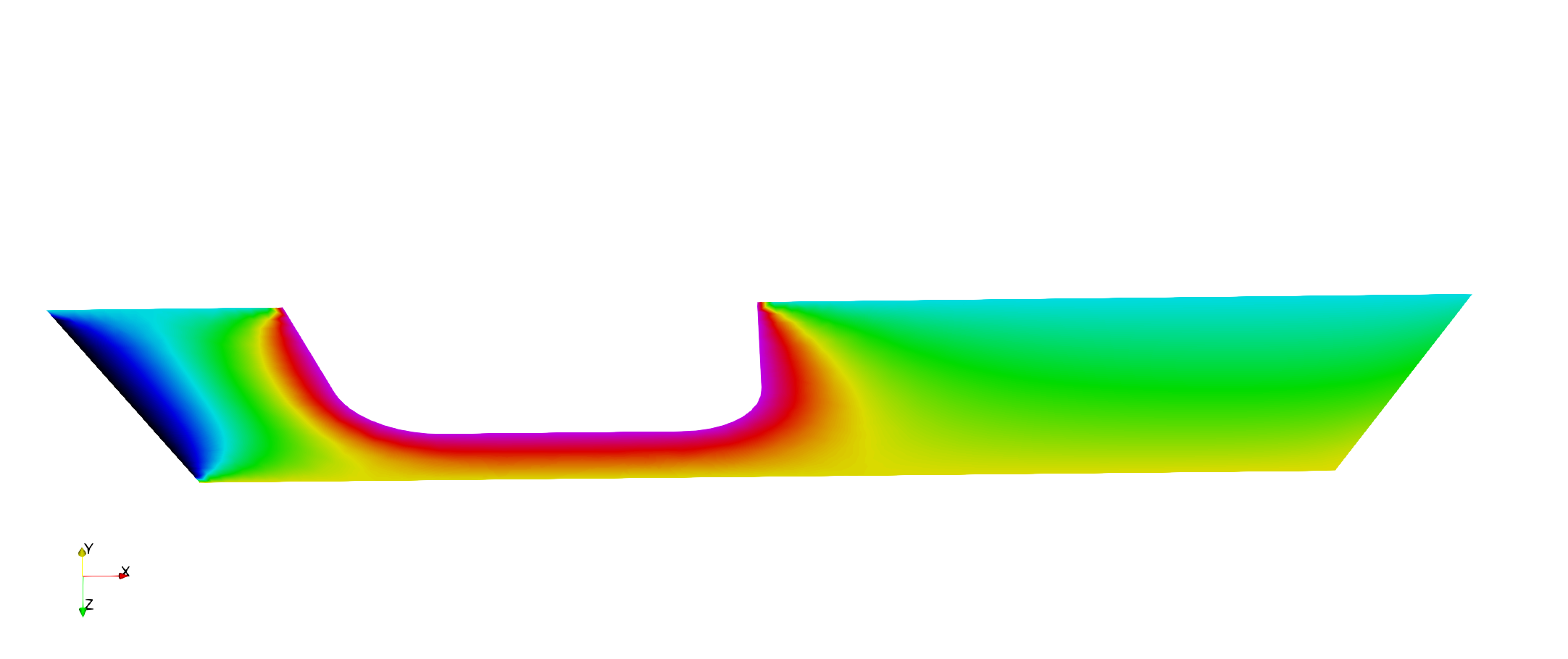}}\\

\subfloat[Pressure distribution at $t=12\, \mu$s \label{fig:valve3P12}]{\includegraphics{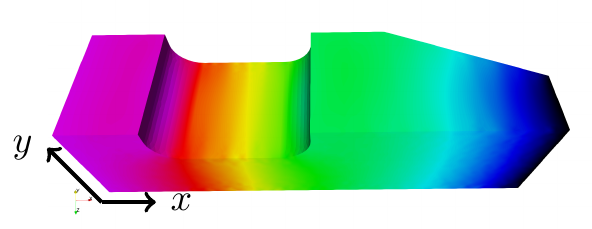}}
\subfloat[Velocity field at $t=12\, \mu$s \label{fig:valve3V12}]{\includegraphics[width=\myW,trim={0cm 0cm 0cm 0cm},clip]{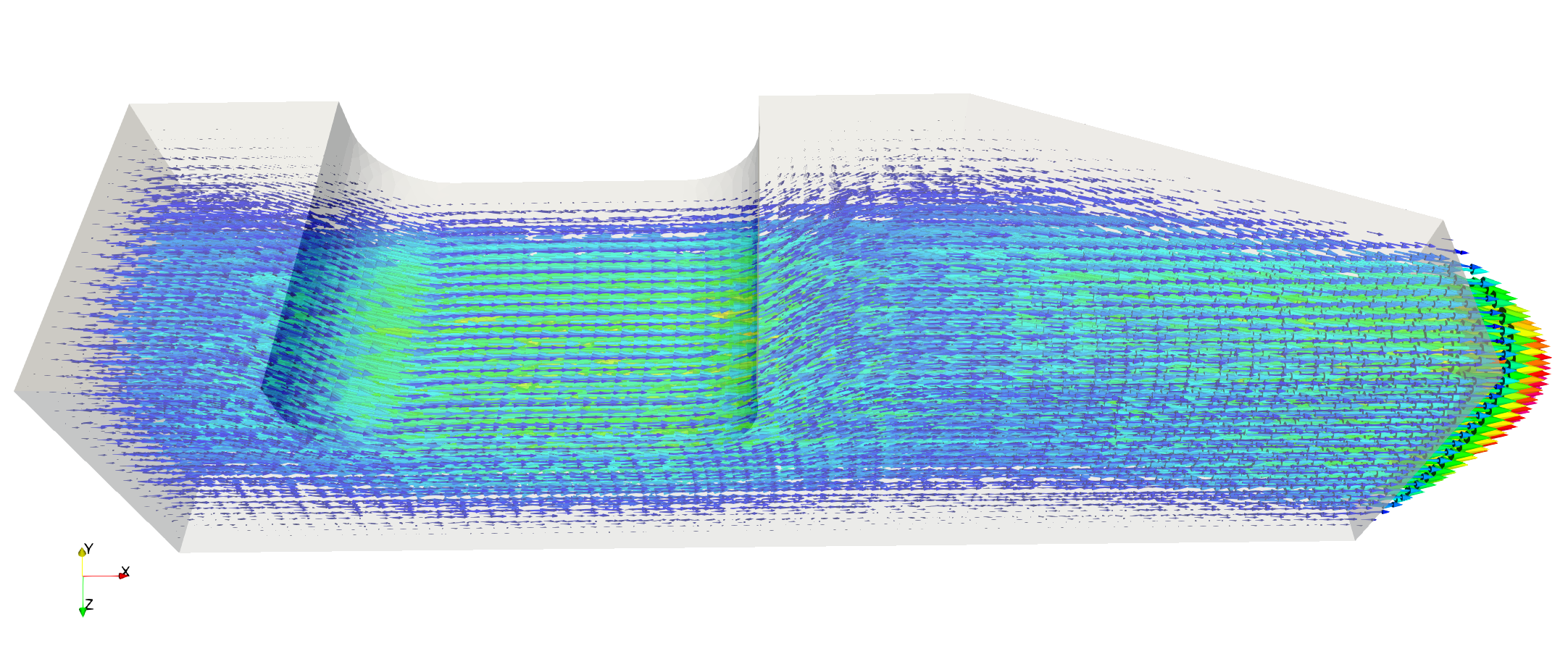}}
\subfloat[Temperature at $z = 3 \, \mu$m, $t=12\, \mu$s \label{fig:valve3T12}]{\includegraphics[width=\myW,trim={0cm 0cm 0cm 0cm},clip]{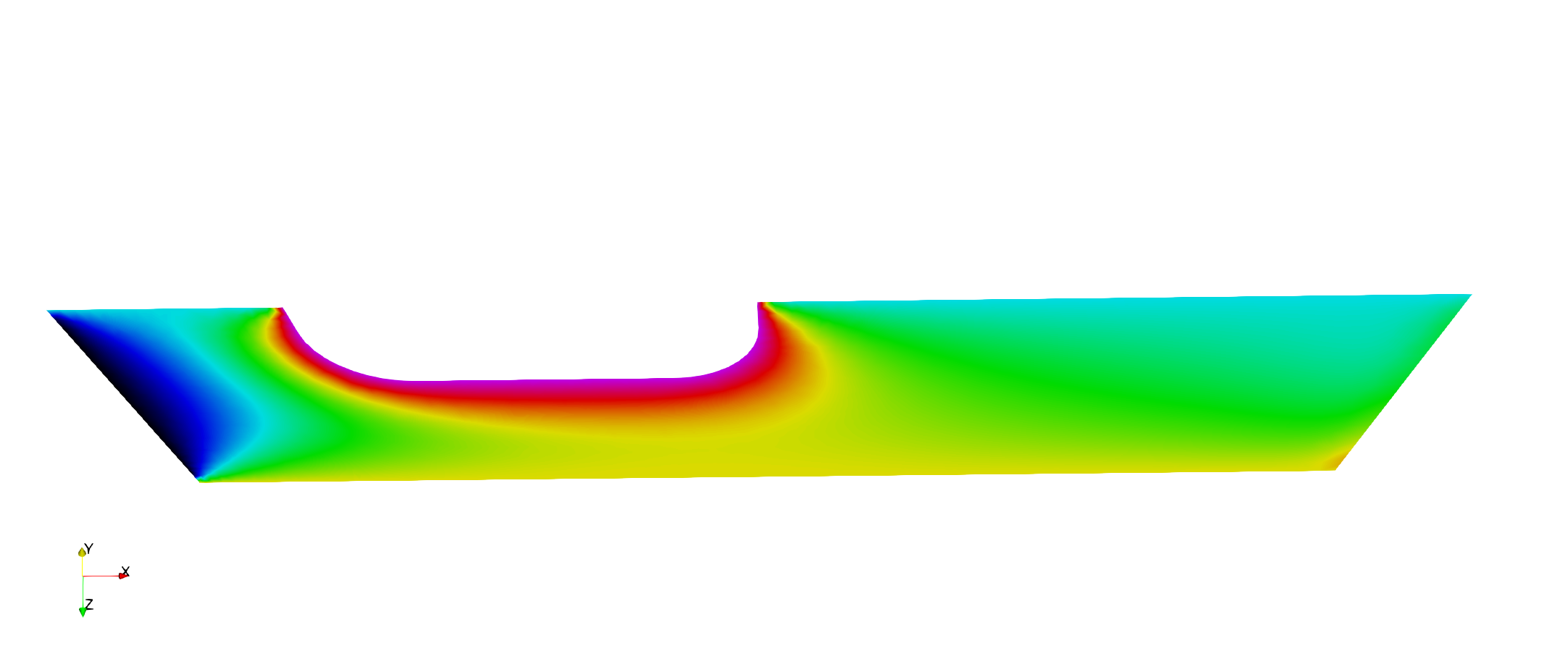}}\\
\caption{Pressure, velocity, and temperature in the closed, half open, and fully opened valve.}
\label{fig:valve4D}
\end{figure}

\subsection{Flow Inspired by Clamped Artery}
\label{ssec:artery}
\begin{figure}
\centering

\subfloat[Side view at $t=0.0\,$s]{\includegraphics{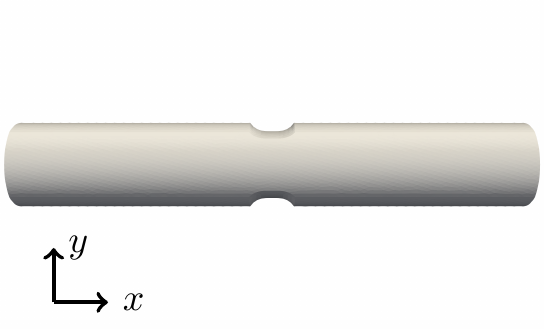}}
\subfloat[Top view $t=0.0\,$s]{\includegraphics{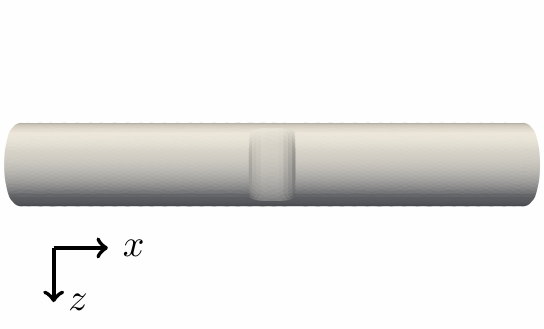}}
\subfloat[Front view at $t=0.0\,$s]{\includegraphics{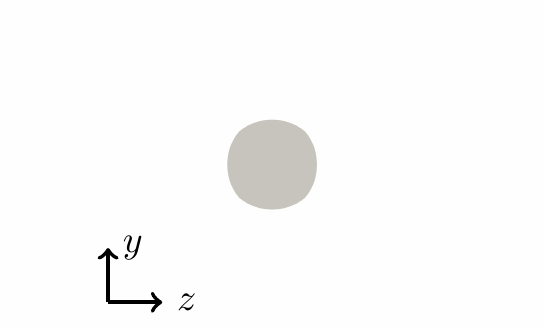}}\\
\subfloat[Side view at $t=0.5\,$s]{\includegraphics{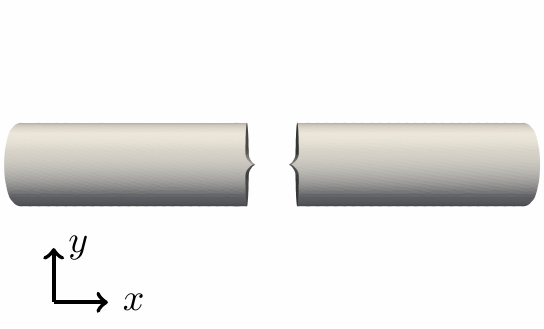}}
\subfloat[Top view at $t=0.5\,$s]{\includegraphics{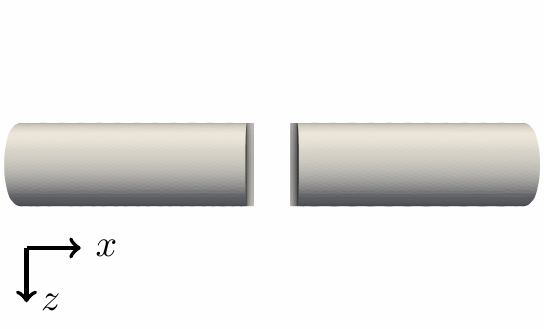}}
\subfloat[Front view at $t=0.5\,$s]{\includegraphics{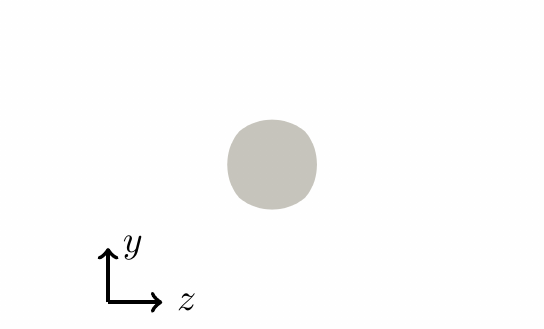}}
\caption{Geometry of clamped artery test case.}
\label{fig:arteryGeo}
\end{figure}

As further application example, we simulate a transient flow through an artery which is temporarily sealed by a clamp and reopened. The considered artery section is \SI{6}{\centi\meter} long and has an approximately circular cross-section with a diameter of \SI{1}{\centi\meter}. The clamp centre is located at $x=\SI{0}{cm}$. Regarding the fluid domain, we assume that the artery volume is displaced by the clamp and returns to the initial approximately circular shape as the clamp is removed. See Figure~\ref{fig:arteryGeo} for several views on the problem geometry in initial and clamped configuration. As time frame for the closing and opening, we choose the duration of one cardiac cycle approximated by \SI{1}{\second}. For the first \SI{0.2}{\second}, the artery is in its initial shape. Over the next \SI{0.2}{\second}, the clamp is applied and seals the artery from \SI{0.4}{\second} until \SI{0.6}{\second}. From \SI{0.6}{\second} until \SI{0.8}{\second}, the artery is reopened and and for the last \SI{0.2}{\second} it is again in its initial shape.

The flow through the deforming domain is driven by a pressure gradient of 40 mmHg, which corresponds to the difference between the minimum and maximum aortic pressure during a cardiac cycle. On the left open boundary, a pressure $ p_{in} = 1 \,\text{atm} + 120 \,\text{mmHg} \approx \SI{1.173e5}{\pascal}$ is prescribed. On the right open boundary, a pressure $p_{out} = 1\, \text{atm} + 80 \,\text{mmHg} \approx \SI{1.12e5}{\pascal}$ is prescribed. We assume that the fluid enters the domain at normal body temperature $T_{in} = \SI{310.15}{\kelvin}$. The colder clamp leads to a temperature variation of  \SI{3}{\kelvin} along the artery wall. The wall temperature is denoted by $T_{W}$. A summary of the corresponding boundary conditions is given in Table~\ref{tab:arteryBc}. Regarding the velocity degrees of freedom, we strongly enforce the no-slip boundary condition on the arterial wall and set the tangential velocity components on the open boundaries to zero. The constant fluid viscosity in this test case, $\nu = 0.04$ poise = \SI{4e-3}{\pascal \second}, is inspired by blood. However, the gas constant, ratio of specific heats and Prandtl number are chosen as in the previous example, leading to a fictitious fluid with a density of roughly \SI[per-mode=symbol]{1}{\kilo\gram\per\cubic\meter}. More accurate models  consider blood as an isothermal generalised Newtonian fluid~\cite{pacheco2020continuous} or even include shear-thinning, viscoelasticity and thixotropy~\cite{owens2006new}. 

\begin{table}
\centering
\caption{Clamped artery test case. Pressure and temperature boundary conditions, viscosity.}
\begin{tabular}{c c c c c c}
\toprule
 $T_{in} \, [\SI{}{\kelvin}]$ & $T_{out} \, [\SI{}{\kelvin}]$ & $T_{W} \, [\SI{}{\kelvin}]$ & $p_{in} \, [\SI{}{\pascal}]$  & $p_{out} \, [\SI{}{\pascal}]$ & $\nu \, [\SI{}{\pascal \second}]$ \\
  \midrule
310.15 & 310.15 &$ 307.15 + 3 \cdot \frac{|x|}{\SI{3}{\centi\meter}}$ & \SI{1.173e5}{}  & \SI{1.12e5}{} & \SI{4e-3}{}\\
 \bottomrule
\label{tab:arteryBc}
\end{tabular}
\end{table}

\subsubsection{Steady Simulation} \label{sssec:steady}

\begin{figure}
\subfloat[Velocity distribution in circular pipe\label{fig:poiseuilleVel}]{
\includegraphics[width=0.53\textwidth,trim={0cm 5cm 10cm 10cm},clip]{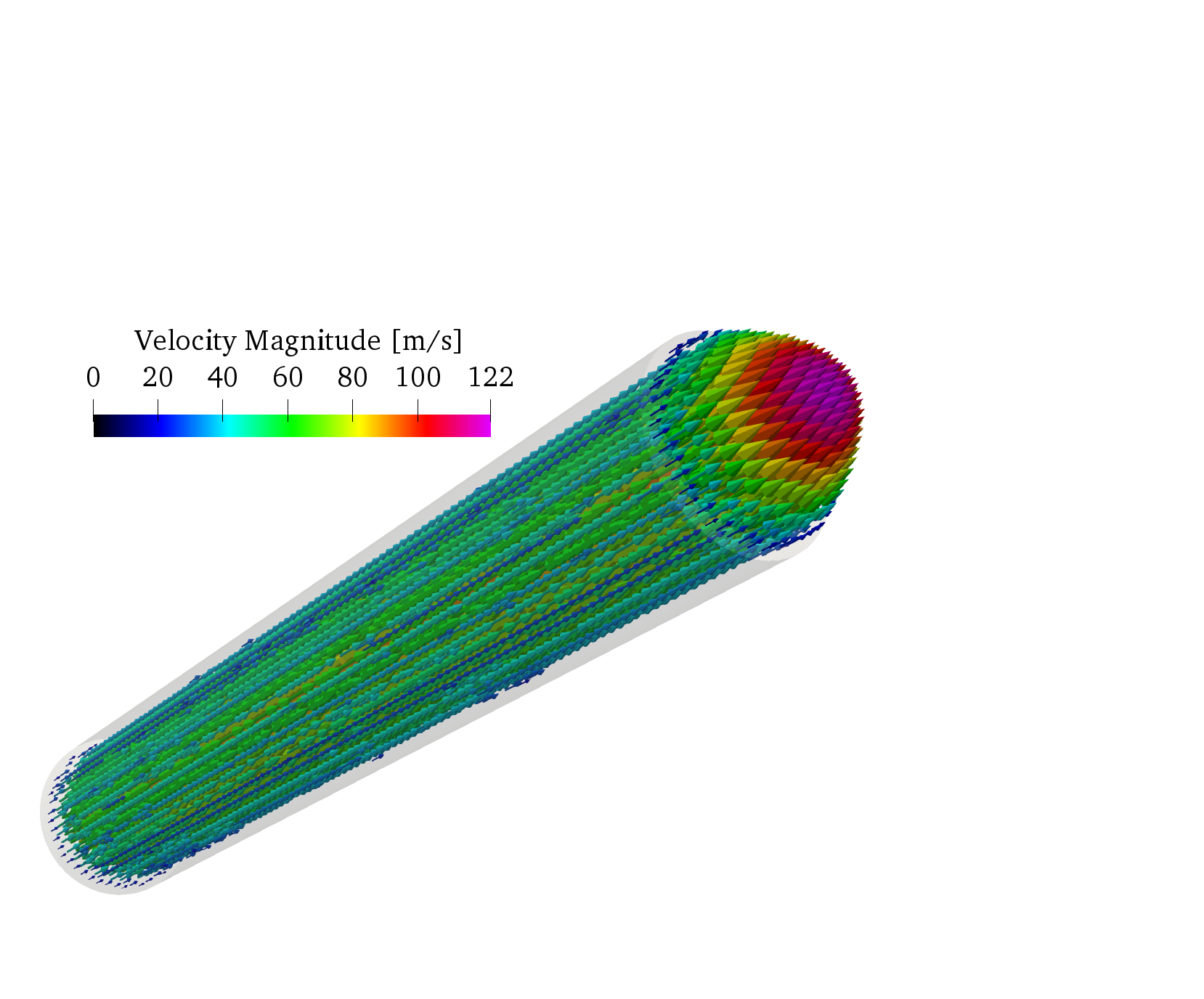}
}
\subfloat[Velocity in axial direction at $x=\SI{2}{\centi \meter}$, $z = \SI{0}{\centi \meter}$ \label{fig:poiseuilleProfile}]{
\includegraphics{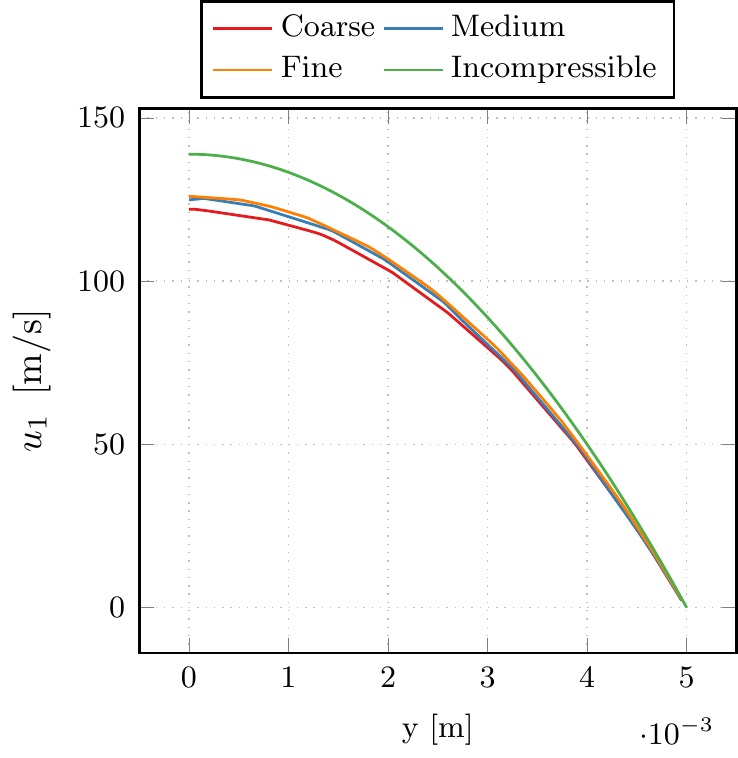}
}
\caption{Compressible Poiseuille flow in circular pipe.}
\label{fig:arteryPoiseuille}
\end{figure}

To validate our formulation for these settings, we perform -- in a first step -- a flow simulation through a straight circular duct. The results are collected in Figure~\ref{fig:arteryPoiseuille}. The velocity distribution is visualised with glyphs showing a paraboloid in qualitative agreement with the Poiseuille paraboloid obtained for an incompressible Hagen-Poiseuille flow. Note that neither on the inflow nor on the outflow boundary the velocity component in normal direction is prescribed. For the pressure driven incompressible flow through a straight circular pipe, there is the well-known analytical solution for the velocity component in axial direction $u_a$~\cite{white2006viscous}, which reads 
\begin{equation}
u_a = -\frac{1}{4\mu} \ddx{p}{x} \left(r_0^2 - r^2 \right) = \SI{138.9}{\meter \per \second} \left(1 - \left(\frac{y}{\SI{5}{\milli \meter}}\right)^2 \right)
\end{equation}
 for the considered configuration. Figure~\ref{fig:poiseuilleProfile} shows that the compressible flow solutions are flatter than $u_a$. The simulation is performed on a coarse, medium, and fine mesh, with \num{111744}, \num{206568}, and \num{652090} tetrahedral elements, respectively. The slight variation of the centerline axial velocity---\SI{122.03}{\meter\per\second} on the coarse mesh (3.1\% smaller than fine), \SI{124.91}{\meter\per\second} on the medium mesh (0.8\% smaller than fine), and \SI{125.98}{\meter\per\second} on the fine mesh---indicates that the coarse mesh resolution is sufficient to obtain a solution within the range of engineering accuracy. 
 
\subsubsection{Transient Simulation} \label{sssec:transient}

\begin{figure}
\centering
\includegraphics[width=0.3\textwidth,trim={0cm 0cm 0cm 0cm},clip]{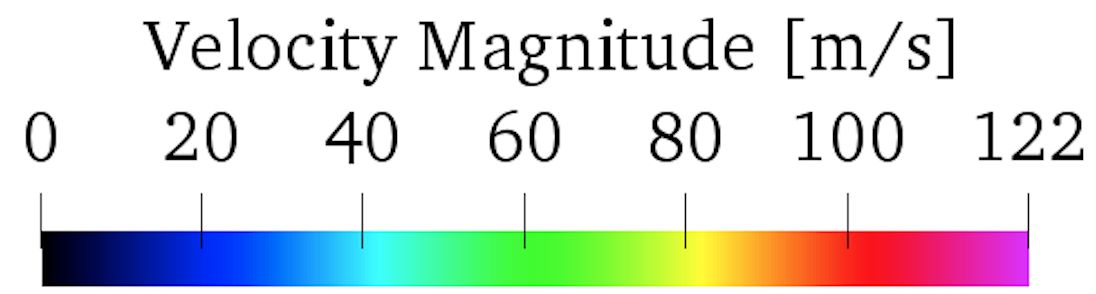}\\
\subfloat[Circular pipe\label{fig:circularVel}]{
\includegraphics[width=0.3\textwidth,trim={2cm 2cm 2cm 2cm},clip]{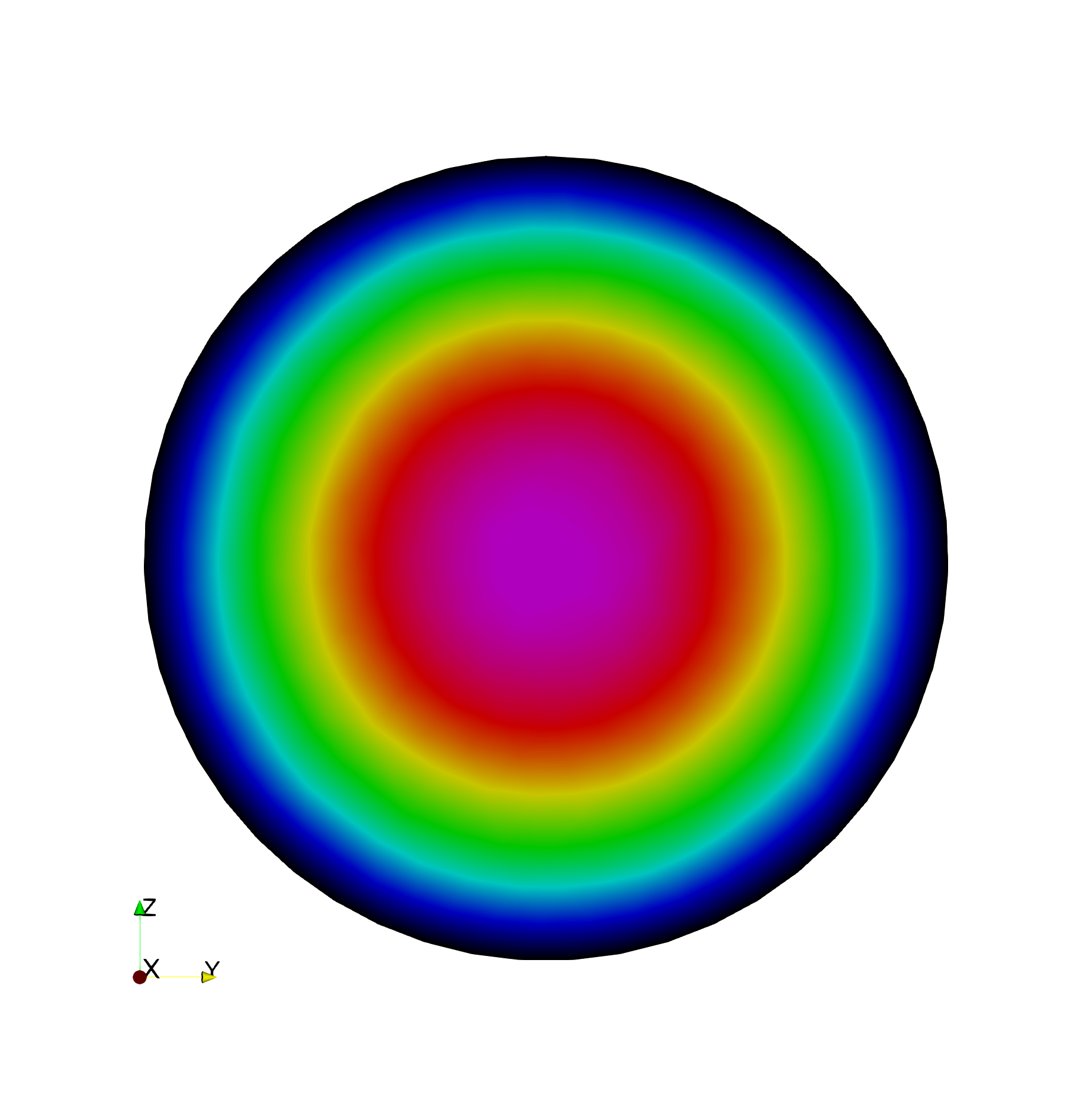} 
}
\subfloat[Approximately circular pipe\label{fig:perturbedVel}]{
\includegraphics[width=0.3\textwidth,trim={2cm 2cm 2cm 2cm},clip]{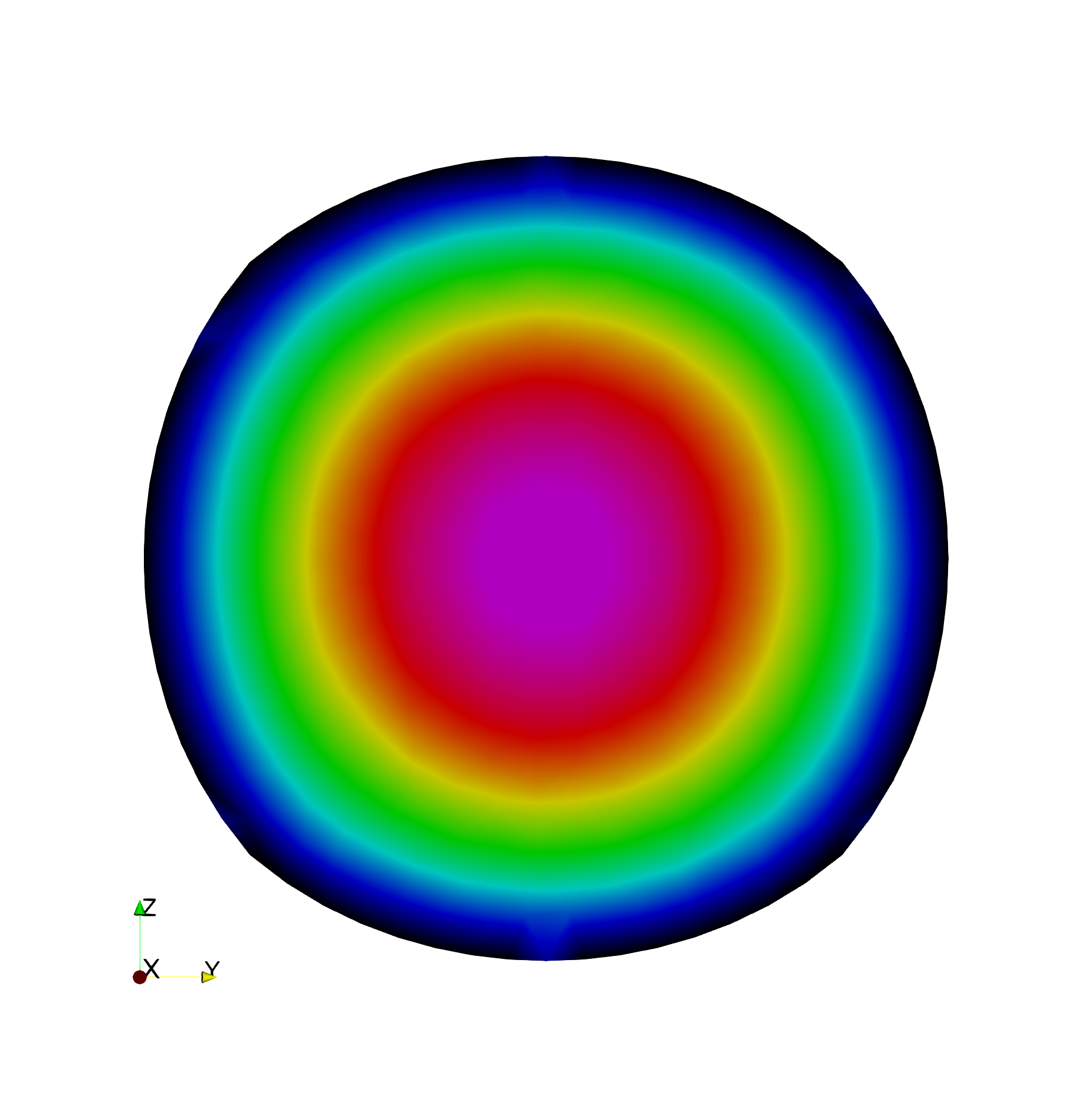} 
}
\caption{Comparison of velocity distribution at $x = \SI{3}{cm}$ between approximated and circular pipe.}
\label{fig:arteryPoiseuilleTrans}
\end{figure}

As second step, we perform a transient simulation on a straight duct with approximately circular cross-section. Starting point for the SST mesh generation is an unstructured tetrahedral $x$-$y$-$t$-Mesh. The mesh has a temporal resolution comparable to the one shown in Figure~\ref{fig:arteryMesh}, but the influence of the clamp is excluded for now. Next, the tetrahedral mesh is extruded in $z$-direction, such that the cross-section in the $y$-$z$-plane forms a square. In $z$-direction, 14 nodes are added during the extrusion, such that the spatial resolution of the resulting mesh is comparable to the coarse pipe mesh discussed above. To obtain the approximately circular cross-section, we perform the 4DEMUM with the space-time coordinates $(x,y,z,t)$ identified with $(x_1, x_2, x_3, x_4)$. Further, the extruded mesh is shifted and scaled, such that $x_1 \in  [-6, 6], \,  x_2 \in [-1,1], \, x_3 \in [-1,1], \, x_4 \in [-1,1]$. On the pentatope mesh boundary $\partial Q_\#^D$, we prescribe the displacements
\begin{equation}
 	\label{eq:artDisp}
 	   \bd_D =  0.9 \left[  \begin{array}{c}
	     0 \\
	     \left( \sqrt{1-\frac{x_3^2}{2}} -1\right) x_2\\ 
	     \left( \sqrt{1-\frac{x_2^2}{2}} -1\right) x_3\\ 
	     0  \end{array} \right],  \quad
\end{equation}
which map the square cross-section in the $x_2$-$x_3$-plane on an approximately circular shape (see Figure~\ref{fig:perturbedVel}). In a final step before the transient simulation, the mesh is shifted and rescaled such that $x \in  [-\SI{3}{\centi\meter}, \SI{3}{\centi\meter}], \,  y \in [-\SI{0.5}{\centi\meter},\SI{0.5}{\centi\meter}], \, z \in [-\SI{0.5}{\centi\meter},\SI{0.5}{\centi\meter}], \, t \in [ \SI{0}{\second} ,\SI{1}{\second}]$. The shifts are performed to allow for relatively simple expressions in the boundary conditions of 4DEMUM, and at the same time allow the fluid simulation to start at $t=0$.

\begin{figure}
\centering
\subfloat[Prescribed inflow pressure\label{fig:presEvo}]{\includegraphics{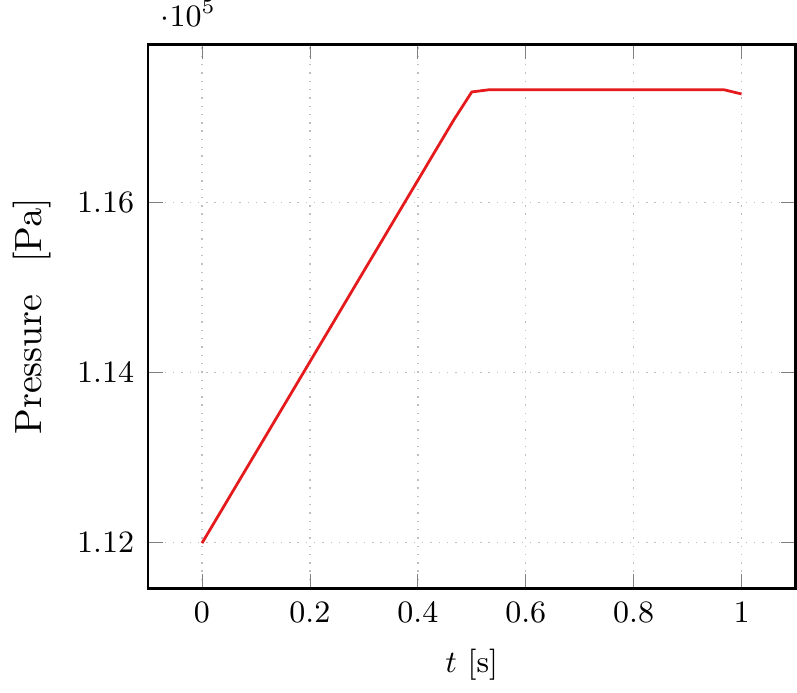}}
\subfloat[Velocity evolution at $x=\SI{-3}{\centi \meter}$, $y = \SI{0}{\centi \meter}$, $z = \SI{0}{\centi \meter}$ \label{fig:veloEvo}]{\includegraphics{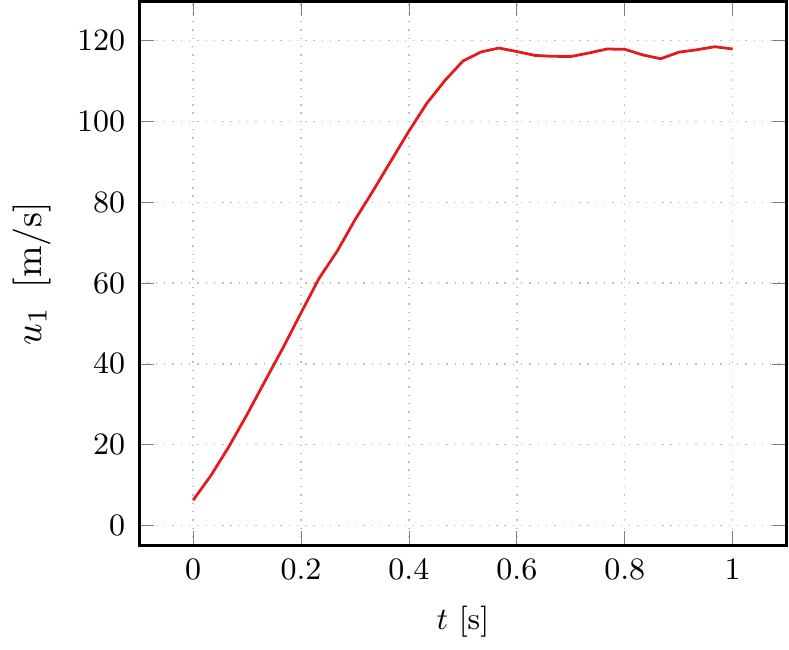}}
\caption{Transient pressure-driven compressible flow in approximately circular pipe.}
\label{fig:arteryTrans}
\end{figure}

In this second simulation, a transient feature is introduced by ramping up the inflow pressure from the initial value $p_{out}$ to $p_{in}$ over the first \SI{0.5}{\second}. For the second \SI{0.5}{\second}, the pressure value is kept constant as shown in Figure~\ref{fig:presEvo}. The flow velocity at the center of the inflow (Figure~\ref{fig:veloEvo}) closely follows the temporal evolution of the pressure value. Note that the computed flow velocity is slightly larger than zero at $t=0$. This is in line with our numerical formulation, which weakly enforces the initial condition (compare Equation~\eqref{eq:weak}). However, the transient nature of the problem is properly captured in the computed flow field. 

With this test case, we also want to explain our choice of an approximately circular cross-section. A perfectly circular cross-section introduces an expected complication in 4DEMUM, i.e., the elements formerly in the corners of the square cross-section attain very large dihedral angles and eventually lead to a tangled mesh. We circumvent this problem by introducing the prefactor in Equation~\eqref{eq:artDisp} to obtain a mesh with an approximately circular cross section in the $y$-$z$-plane. 

\begin{figure}
\centering
\includegraphics{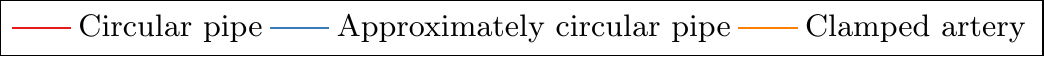}\\
\subfloat[Velocity $u_1$ at $x=\SI{2}{\centi \meter}$, $z = \SI{0}{\centi \meter}$, $t = \SI{0.9}{\second}$ \label{fig:arteryCompProfile}]{
\includegraphics{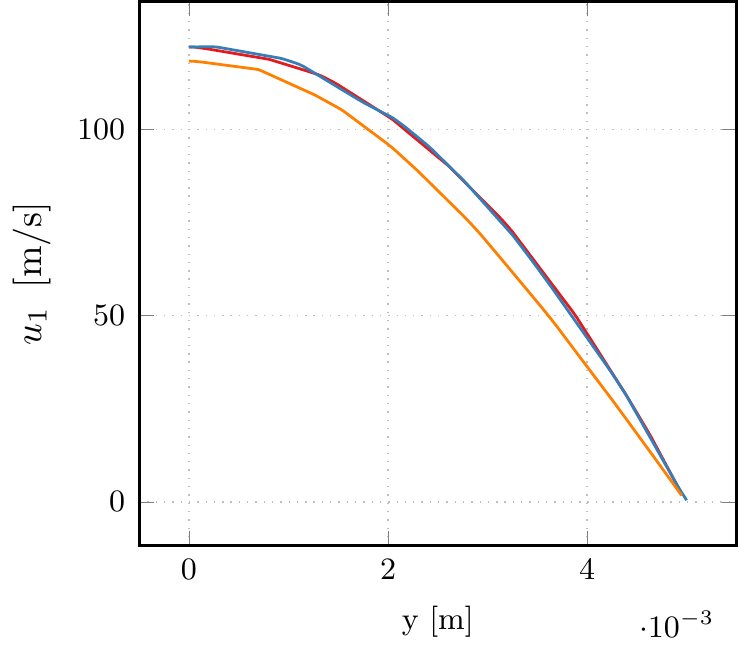}
} \quad
\subfloat[Pressure at centre line $y=\SI{0}{\centi \meter}$, $z = \SI{0}{\centi \meter}$, $t = \SI{0.9}{\second}$ \label{fig:arteryCompPressure}]{
\includegraphics{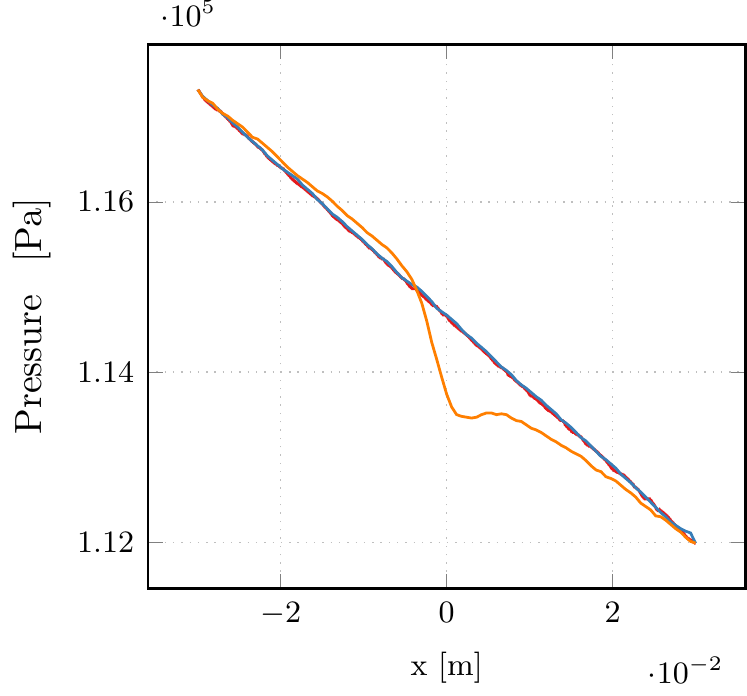}
}
\caption{Comparison between circular pipe, approximately circular pipe, and artery geometry.}
\label{fig:arteryComp}
\end{figure}

The ``missing 10\%`` towards the perfect circular cross-section have hardly any influence on the flow field as shown in Figure~\ref{fig:arteryPoiseuilleTrans}. A more quantitative comparison is presented in Figure~\ref{fig:arteryComp}. The parabolic velocity profile in radial direction (Figure~\ref{fig:arteryCompProfile}) as well as the linear pressure decay along the pipe axis (Figure~\ref{fig:arteryCompPressure}) are obtained independently of the approximation of the circular cross-section. Therefore, we consider an approximately circular cross-section of the artery in the following. 

\subsubsection{Transient Simulation with Topology Change} \label{sssec:topology}

\begin{figure}
\centering
\includegraphics{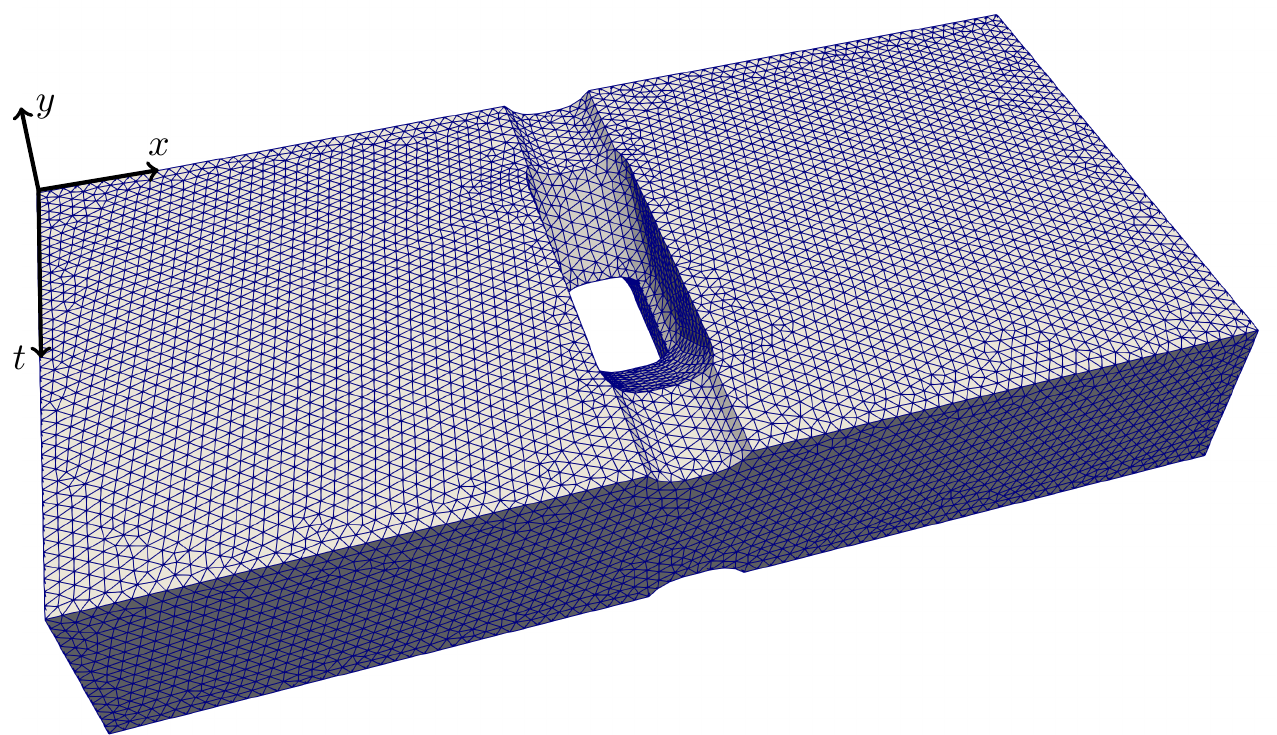}
\caption{Clamped artery test case. $x$-$y$-$t$-Mesh}
\label{fig:arteryMesh}
\end{figure}

As third and final step, we consider the transient simulation with topology change of the spatial computational domain. We now include the topology change caused by the clamp in the $x$-$y$-$t$-mesh shown in Figure~\ref{fig:arteryMesh}. The further pentatope mesh generation steps of extrusion, connectivity generation, and elastic deformation are performed as in the previous example; the boundary displacements are given in Equation~\eqref{eq:artDisp}.

The subsequent finite element flow simulation was performed in on 240 cores using a distributed memory parallelisation based on MPI. It took 16 minutes of wall clock time. Figure~\ref{fig:artery} presents the simulation results at $t =$ \SI{0.5}{\second}, \SI{0.7}{\second}, and \SI{0.9}{\second}. For the fully clamped case, we obtain two separate domains and negligible flow velocities (Figure~\ref{fig:art3P5}). In the absence of flow, the temperature distribution in the fluid is an interpolation of the temperature prescribed on the domain boundaries (Figure~\ref{fig:art3T5}). Note that on the right most boundary ($x=\SI{3}{\centi \meter}$) no temperature boundary condition is prescribed because this part is an outflow boundary during most of the simulation. However, from $t=\SI{0.5}{\second}$ until $t=\SI{0.7}{\second}$, back-flow across this boundary introduces a small disturbance in the temperature field. When the artery is reopened to roughly half of the total diameter, a strong pressure gradient across the clamp region accelerates the flow in this area (Figure~\ref{fig:art3P7}). Further at $t=\SI{0.7}{\second}$, the temperature distribution is strongly influenced by the flow field as well as the cooler clamp (Figure~\ref{fig:art3T7}). In the open configuration (Figure~\ref{fig:art3P9} and ~\ref{fig:art3T9}), we observe a linear pressure decrease from the inflow to the outflow and a parabolic velocity profile everywhere except for the clamp region. The comparison to the velocity profile of the approximately circular pipe shows that the flow speed computed for the clamped artery is slower at $x=\SI{2}{\centi \meter}$, $z=\SI{0}{\centi \meter}$, $t=\SI{0.9}{\second}$ (Figure~\ref{fig:arteryCompProfile}). This is to be expected as the same pressure gradient has to overcome the additional obstacle of the clamp region (Figure~\ref{fig:arteryCompPressure}).

\renewcommand{\myW}{0.47\textwidth}
\begin{figure}
\centering
\includegraphics[width=0.35\textwidth,trim={23.5cm 26cm 22.5cm 10cm},clip]{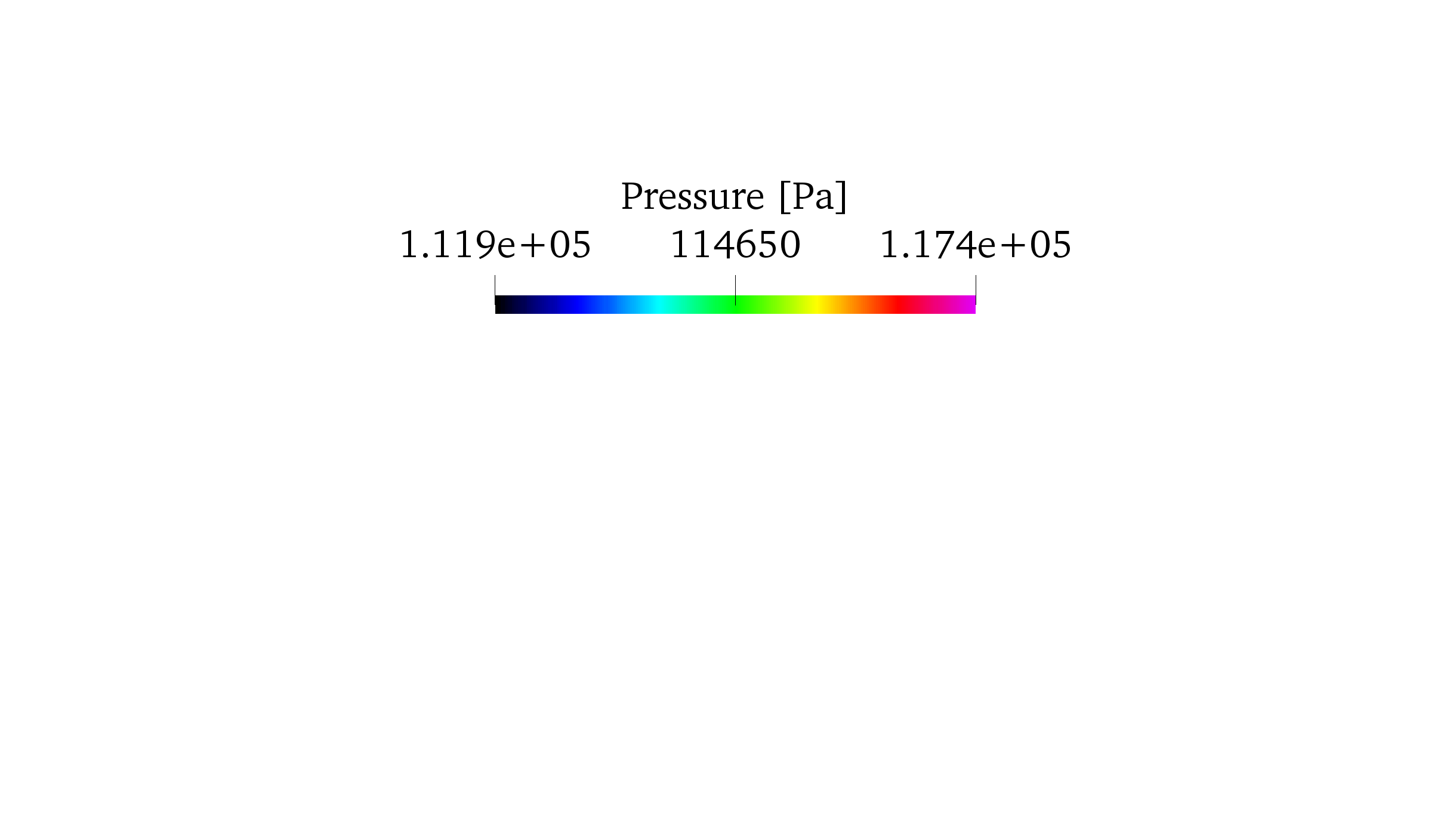}\qquad
\includegraphics[width=0.35\textwidth,trim={23.5cm 26cm 22.5cm 10cm},clip]{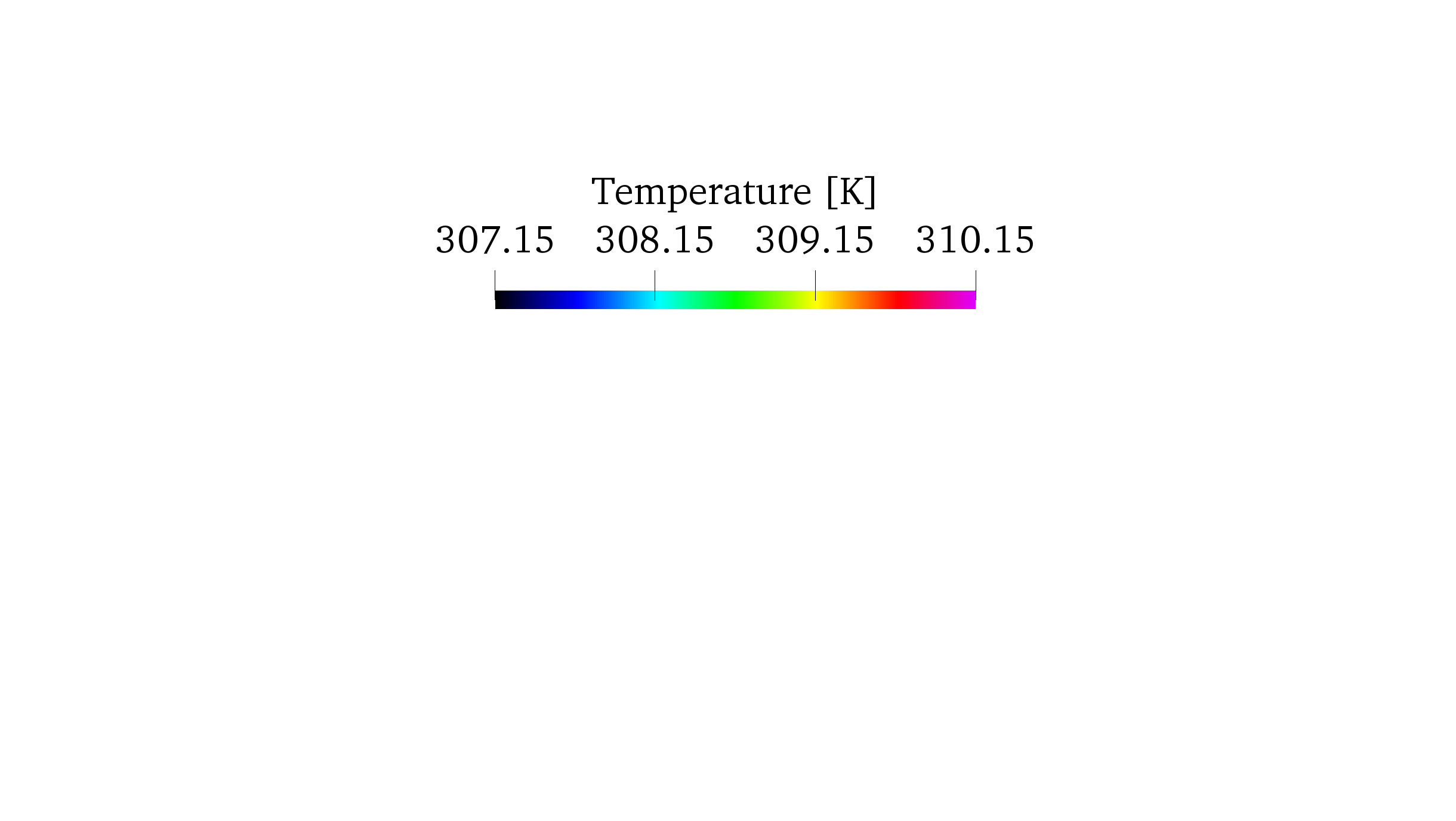}

\subfloat[Pressure distribution on velocity glyphs at $t=0.5\,$s \label{fig:art3P5}]{\includegraphics[width=\myW,trim={0cm 3cm 0cm 0cm},clip]{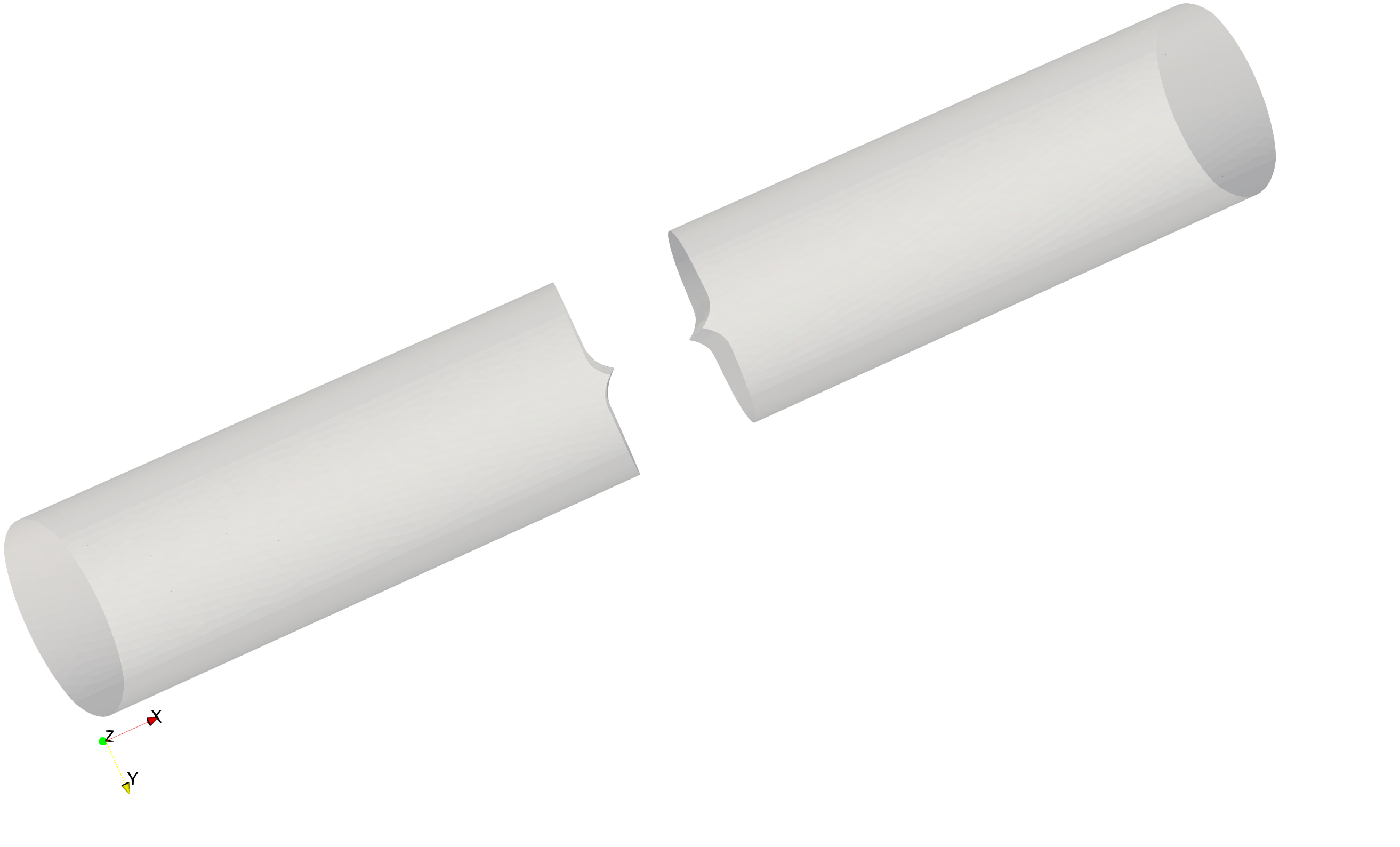}}\quad
\subfloat[Temperature at $z=0$, $t=0.5\,$s \label{fig:art3T5}]{\includegraphics[width=\myW,trim={0cm 3cm 0cm 0cm},clip]{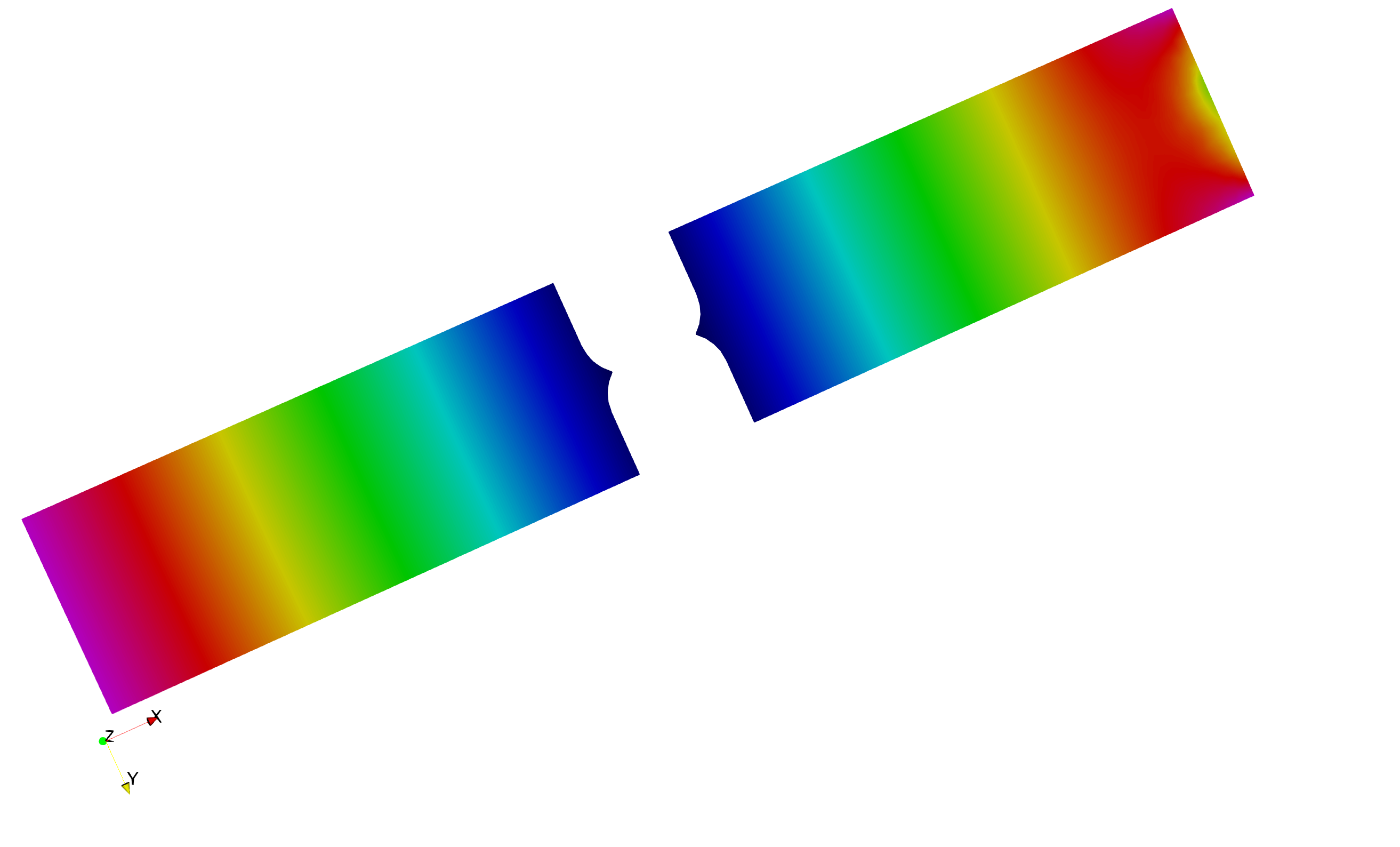}}

\subfloat[Pressure distribution on velocity glyphs at $t=0.7\,$s \label{fig:art3P7}]{\includegraphics[width=\myW,trim={0cm 3cm 0cm 0cm},clip]{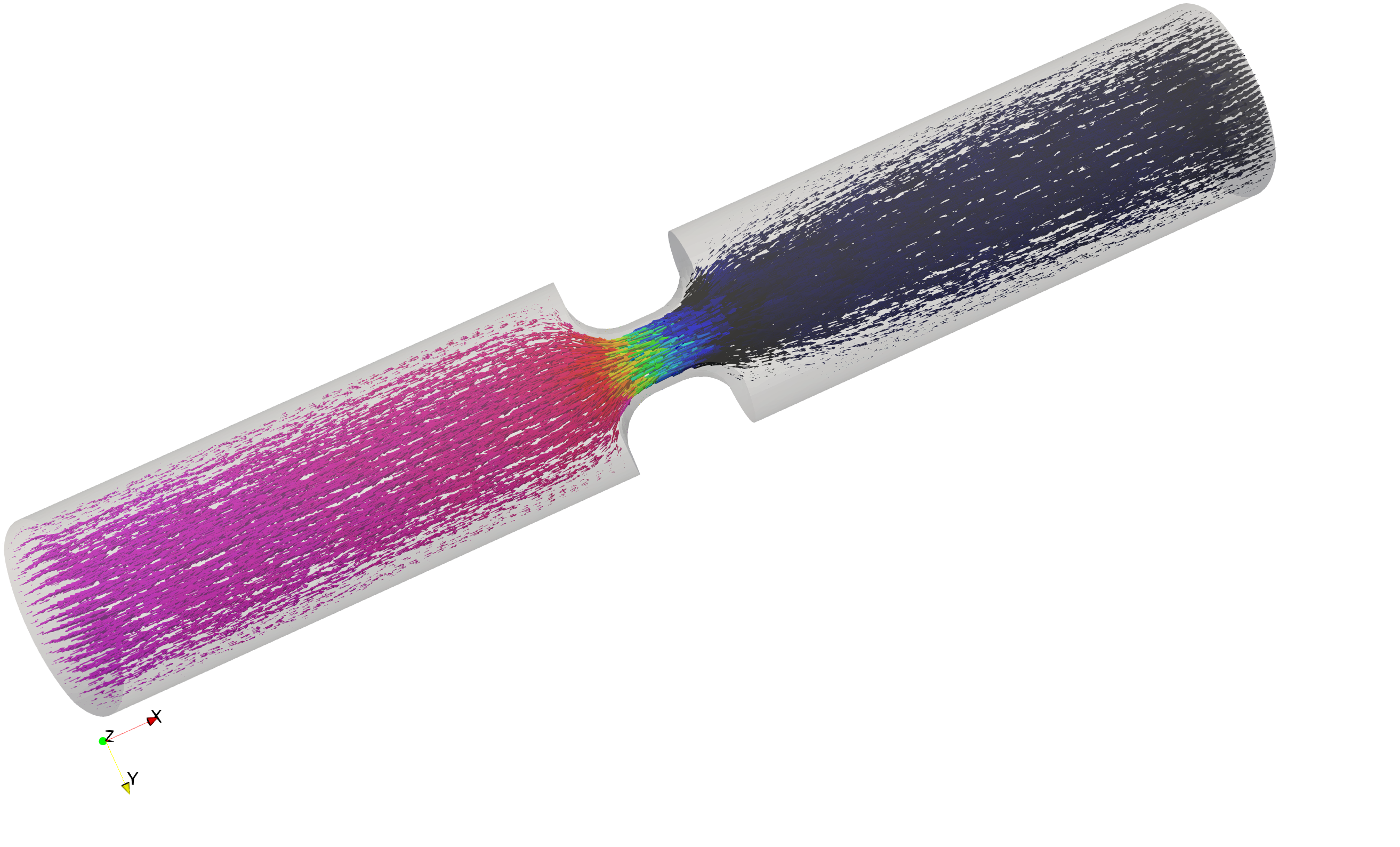}}\quad
\subfloat[Temperature at $z=0$, $t=0.7\,$s \label{fig:art3T7}]{\includegraphics[width=\myW,trim={0cm 3cm 0cm 0cm},clip]{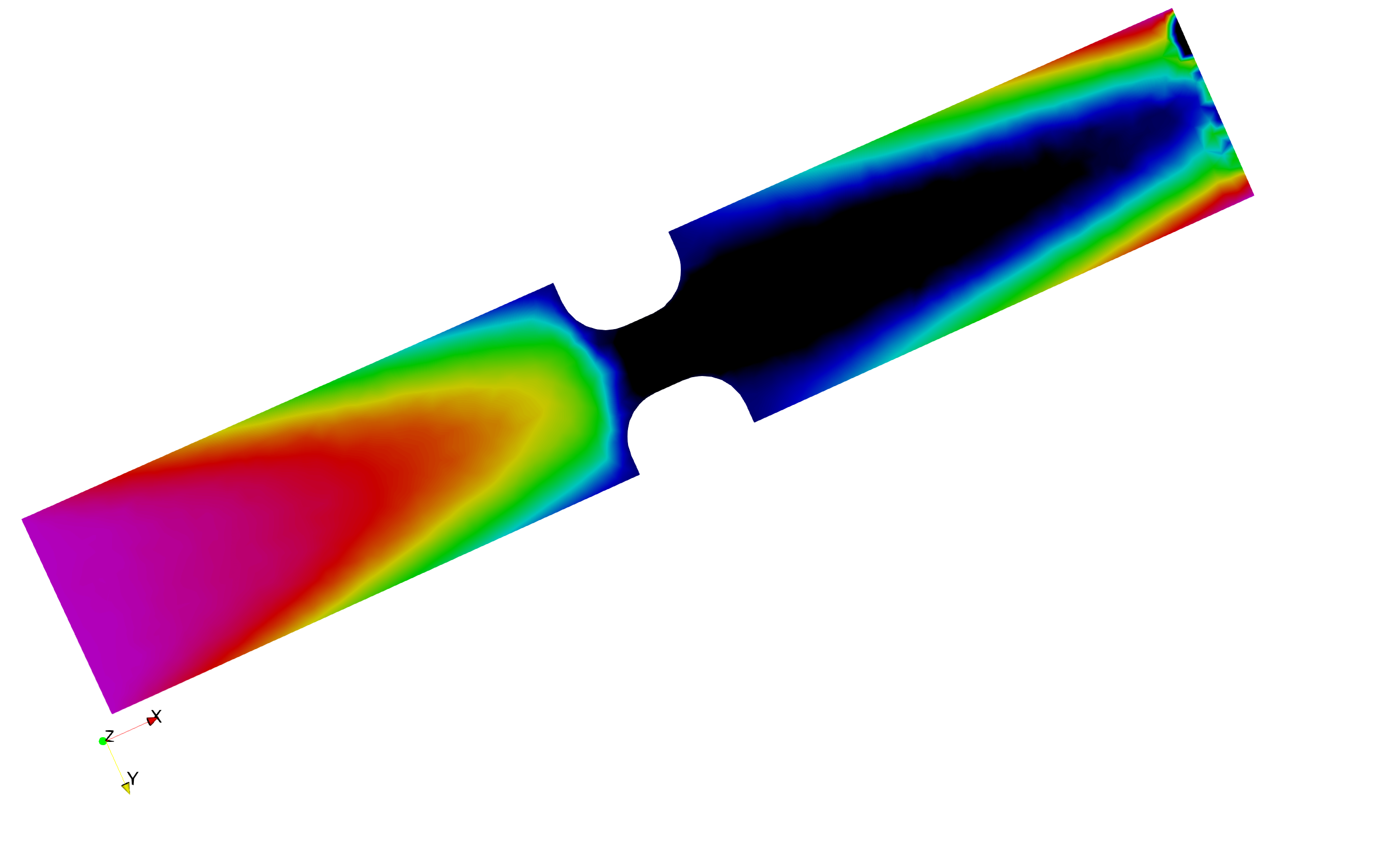}}

\subfloat[Pressure distribution on velocity glyphs at $t=0.9\,$s \label{fig:art3P9}]{\includegraphics[width=\myW,trim={0cm 3cm 0cm 0cm},clip]{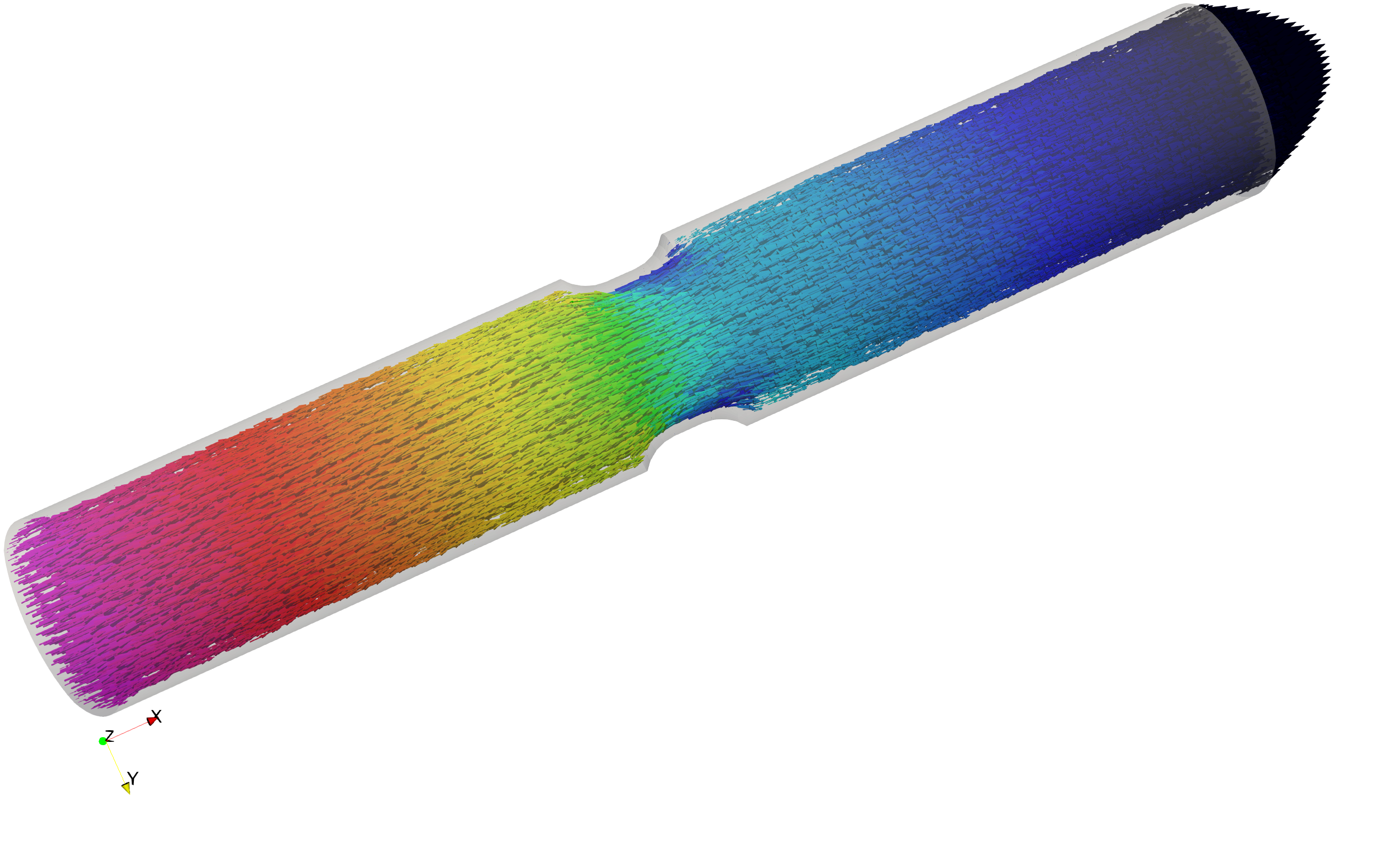}}\quad
\subfloat[Temperature at $z=0$, $t=0.9\,$s \label{fig:art3T9}]{\includegraphics[width=\myW,trim={0cm 3cm 0cm 0cm},clip]{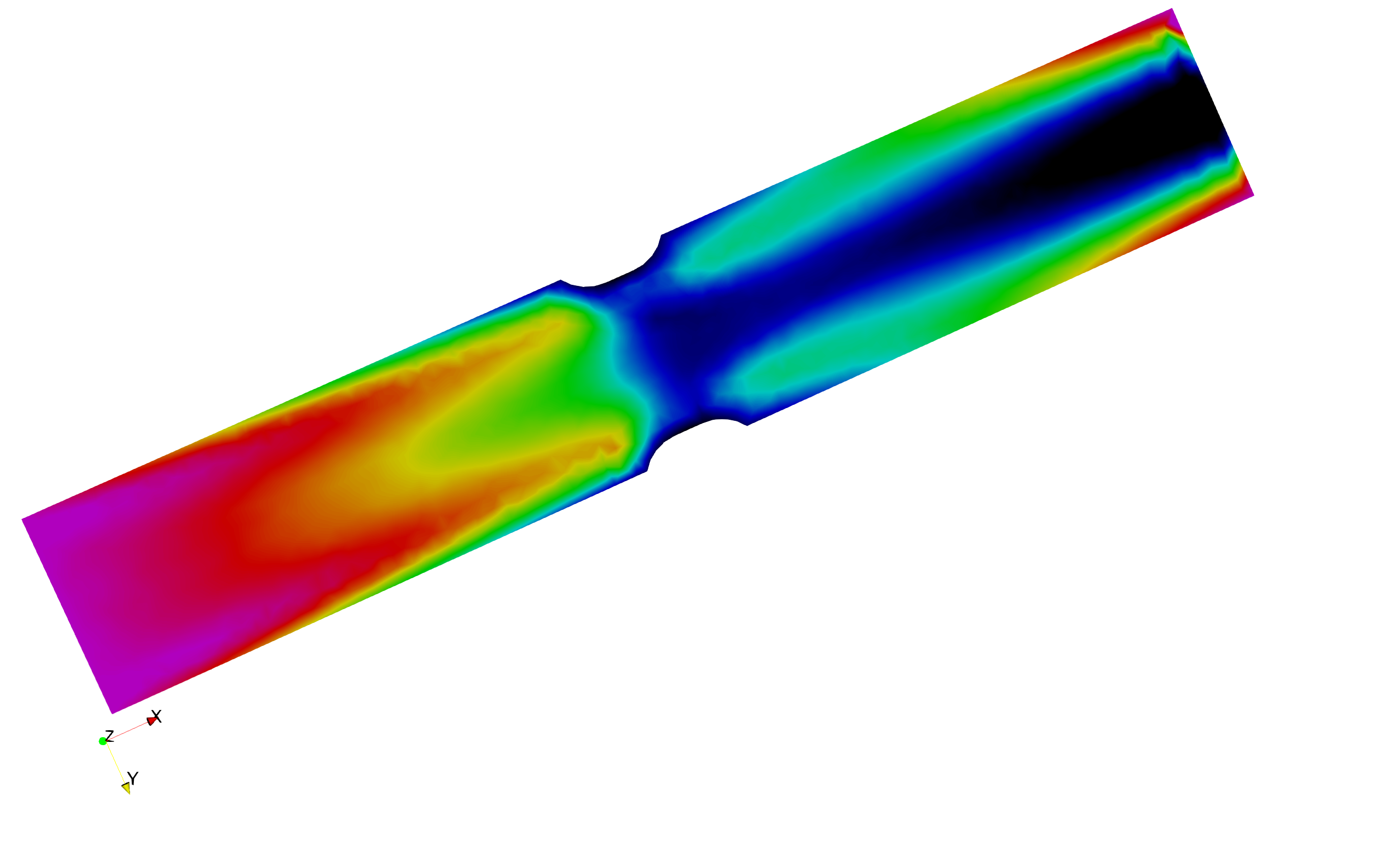}}

\caption{Pressure, velocity, and temperature in the fully clamped, partially, and fully opened artery.}
\label{fig:artery}
\end{figure}

\section{Conclusion and Outlook}
\label{sec:conclusion}
In this paper, we presented the generation and application of four-dimensional SST meshes that allow for a boundary conforming discretisation of spatial domains with time variant topology. To produce pentatope meshes of complex geometries, the elastic mesh update method was extended to four dimensions. We described the integration of 4DEMUM into the simulation workflow. Namely, we described the steps from the space-time geometry to the SST mesh of the physical domain, as well as the post-processing steps to visualise the solution on the four-dimensional mesh as series of data sets on three-dimensional meshes. The workflow was successfully applied to two test cases featuring geometries of a valve and a clamped artery. The transient three-dimensional flow solutions validate the mesh generation in the sense that proper pentatope finite element meshes are obtained.

As an outlook, the enhanced meshing capabilities open up a path to parallel-in-time computations on complex domains. We applied domain-decomposition not only to the spatial domain, but to the complete space-time domain. Another promising application area are FSCI simulations with topology changes on boundary conforming meshes. For this application, the topology changes have to be included in the SST mesh, however, spatial and temporal position of the topology changes can be determined in the course of a coupled FSI simulation and adjusted using 4DEMUM.

\end{document}